\documentclass[final]{siamltex}
\usepackage{graphicx}
\usepackage{psfrag}
\usepackage{amsmath,rotating}
\usepackage{amsfonts,amsbsy,dsfont}
\usepackage{amstext,amssymb,epsf}
\usepackage[dvips]{epsfig}
\usepackage{subfigure}
\usepackage{csquotes}

\def\G{\mathrm{G}}
\def\xx{\mathrm{x}}
\def\vv{\mathrm{v}}
\def\bfA{\mathbf{A}}
\def\bfM{\mathbf{M}}
\def\bfF{\mathbf{F}}
\def\ss{\mathrm{s}}
\def\sd{\, \mathrm{d}}
\def\rd{\mathrm{d}}

\def\app{{a}}

\newcommand{\bv}{{\bar v}}
\newcommand{\bvd}{{\overline{v^2}}}

\title{Multivariate Gaussian extended quadrature method of moments for turbulent disperse multiphase flow}

\author{C.~Chalons\footnotemark[2], R.~O.~Fox\footnotemark[3] \footnotemark[4]~, F.~Laurent\footnotemark[4] \footnotemark[5]~, M.~Massot\footnotemark[4] \footnotemark[5]~, \and A.~Vi\'{e}\footnotemark[4]}

\begin{document}

\maketitle

\renewcommand{\thefootnote}{\fnsymbol{footnote}}
\footnotetext[2]{LMV - UMR 8100, Univ. Versailles Saint-Quentin-en-Yvelines, UFR des Sciences, B\^atiment Fermat, 45 avenue des Etats-Unis, 78035 Versailles cedex, France}
\footnotetext[3]{Department of Chemical and Biological Engineering, Iowa State University, Ames, IA 50011-2230, USA}
\footnotetext[4]{Laboratoire EM2C, UPR CNRS 288, CentraleSup\'elec, Universit\'e Paris-Saclay, Grande Voie des Vignes, 92295 Ch\^atenay Malabry, France}
\footnotetext[5]{F\'ed\'eration de Math\'ematiques de l'Ecole Centrale Paris FR CNRS 3487, CentraleSup\'elec, Grande Voie des Vignes, 92295 Ch\^atenay Malabry, France}
\renewcommand{\thefootnote}{\arabic{footnote}}

\begin{abstract}
The present contribution introduces a fourth-order moment formalism for particle trajectory crossing (PTC) in the framework of multiscale modeling of disperse multiphase flow.  
In our previous work, the ability to treat PTC was examined with direct-numerical simulations (DNS) using either quadrature reconstruction based on a sum of Dirac delta functions denoted as Quadrature-Based Moment Methods (QBMM) in order to capture large scale trajectory crossing, or by using low order hydrodynamics closures in the Levermore hierarchy denoted as Kinetic-Based Moment Methods (KBMM) in order to capture small scale trajectory crossing. 
Whereas KBMM leads to well-posed PDEs and has a hard time capturing large scale trajectory crossing for particles with enough inertia, QBMM based on a discrete reconstruction suffers from singularity formation and requires too many moments in order to capture the effect of PTC at both small scale and large scale both to small-scale turbulence as well as free transport coupled to drag in an Eulerian mesoscale framework. 
The challenge addressed in this work is thus twofold: first, to propose a new generation of method at the interface between QBMM and KBMM with less singular behavior and the associated proper mathematical properties, which is able to capture both small scale and large scale trajectory crossing, and second to limit the number of moments used for applicability in 2-D and 3-D configurations without losing too much accuracy in the representation of spatial fluxes. 
In order to illustrate its numerical properties, the proposed Gaussian extended quadrature method of moments (Gaussian-EQMOM) is applied to solve 1-D and 2-D kinetic equations representing finite-Stokes-number particles in a known turbulent fluid flow.
\end{abstract}




\pagestyle{myheadings}
\thispagestyle{plain}
\markboth{C.~Chalons, R.~O.~Fox, F.~Laurent, M.~Massot \& A. Vi\'e }{MULTIVARIATE GAUSSIAN-EQMOM}

\section{Introduction}\label{s:intro}

The physics of inertial particles in a carrier fluid (e.g., fluidized beds, sprays, alumina particles in rocket boosters) can be described by a number density function (NDF) satisfying a kinetic equation. Solving such a kinetic equation relies on either a sample of discrete numerical parcels through a Lagrangian Monte-Carlo approach or on a moment approach resulting in a Eulerian system of conservation laws on velocity moments. For the latter, the main difficulty for particles with high Knudsen numbers (i.e., weakly collisional flows) where the velocity distribution can be very far from equilibrium, is the closure of the free-transport term in the kinetic equation. 

In the field of high-Knudsen rarefied gases, Grad~\cite{grad1949} introduced a perturbation of the equilibrium distribution using Hermite polynomials whose coefficients are found by transporting additional moments, whereas Levermore introduced a naturally well-posed hierarchy of hyperbolic systems of conservation equations \cite{Lev96}.  More recently, Struchtrup and Torrilhon~\cite{struchtrup2003} derived a regularization of the Grad 13-moment system; however, such systems lose the hyperbolic character of the original kinetic equation (see for instance~\cite{torrilhon2000}). To address this issue, methods based on a realizable presumed NDF have been suggested. Torrilhon~\cite{Tor10} used a Pearson-IV distribution as a basis NDF, which extends the hyperbolic region of the 13-moment method, but still is not globally hyperbolic. A class of hyperbolic methods was suggested by Levermore and Morokoff using entropy maximization (EM)~\cite{Lev96}. For moments up to second order, i.e., multivariate Gaussian distributions, the algorithm is fast and easy to implement, and has been used in rarefied gases~\cite{mcdonald2008,mcdonald2014} as well as in particulate flows~\cite{vie2013cicp,sabat2014jcmf,sabat16,sabatth16}. For higher-order moments and multi-dimensional systems, McDonald and co-workers~\cite{mcdonald2013,mcdonald2013b} propose CFD strategies based on EM, yielding a high level of accuracy but with the drawback of an expensive inversion algorithm; besides, as investigated in \cite{hauck2008}, the structure of the moment space leads to some difficulties. This class of methods is referred to as the class of Kinetic-Based Moment Methods.

Another way to proceed is to use Quadrature-Based Moment Methods (QBMM) where the higher-order moments required for closure are evaluated from the lower-order transported moments using multi-dimensional quadratures. In our previous work, we developed multivariate quadratures in the form of a sum of Dirac delta functions in velocity phase space (see \cite{yuan2011,Kah10a} and the references therein). Such an approximation also allows for a well-behaved kinetic numerical scheme in the spirit of \cite{Sdc09_us,massot09_vki,deChaisemartin2009} inspired from \cite{Bouchut03} where the fluxes in a cell-centered finite-volume formulation are directly evaluated from the knowledge of the quadrature abscissas and weights with guaranteed realizability conditions. Such a quadrature approach has been shown to capture particle trajectory crossing (PTC) in direct-numerical simulations (DNS) where the distribution in the exact kinetic equation remains at all times in the form of a sum of Dirac delta functions \cite{Kah10a,kah10,ctrfreret10,yuan2011,vikas2011}, but lead to weakly hyperbolic system of conservation equation and can develop artificial singularities when the number of crossing goes beyond the number of Dirac delta function allowed as well as difficulties at the boundary of the moment space \cite{kah11_cms}.  

In a turbulent fluid, the effect of the turbulence on the particles through transport and drag leads to dispersion in velocity phase space due to the large number of PTC occurring over a wide range of length and time scales. Nevertheless, classical large scale PTC still occurs for large enough Stokes numbers because the small-scale dispersion is not strong enough to ``randomize'' the particle velocities resulting from free transport and drag is not strong enough to prevent such large scale trajectory crossing. However, capturing both small scale leading to velocity dispersion as well as large scale PTC  requires a large number of quadrature nodes using a delta-function representation. Moreover, such QBMM result in entropic weakly hyperbolic systems of conservation laws and to the formation of $\delta$-shock singularities, the mathematical structure of which is studied in \cite{kah11_cms}. The purpose of the present contribution is to introduce an extended quadrature-based reconstruction of the NDF for the closure of the free-transport term in the kinetic equation and drag, in between QBMM and KBMM methods, which also allows us to naturally account for velocity dispersion and large scale trajectory crossing. The proposed continuous representation of the NDF allows us both to limit the number of unknowns in multi-dimensional configurations, and to regularize the resulting system of conservation equations, while still being able to capture large scale PTC and velocity dispersion.

One-dimensional (1-D) multivariate Gaussian moment methods were introduced in \cite{ctr10b,laurent12brief} as Multi-Gaussian KBMM.  A generalization of these methods is proposed in \cite{yuan2012} and called the extended quadrature method of moments (EQMOM).
Here we analyze what we refer to as Gaussian-EQMOM, extend it to 2- and 3-D phase space, and apply it to the solution of the kinetic equation describing inertial particles in a continuous fluid phase.  The remainder of the work is organized as follows.  In \S\ref{s:omgq} a brief introduction to 1-D Gaussian-EQMOM is provided.  In \S\ref{s:quad-1D}, we provide an in-depth description of the application to 1-D kinetic equations and the mathematical properties of two-node Gaussian-EQMOM.  In \S\ref{s:gquad-2D}, we describe the extension of Gaussian-EQMOM to a 2-D phase space using the conditional quadrature method of moments (CQMOM). 
In \S\ref{s:quad-2D} we describe the application of Gaussian-ECQMOM to 2-D kinetic equations. Example applications are provided in \S\ref{s:examples}. Finally, conclusions are drawn in \S\ref{s:conclusion}.

\section{One-dimensional Gaussian-EQMOM}\label{s:omgq}

Consider an NDF $f(v)$ defined for $v \in \mathbb{S}$ where $\mathbb{S} \subseteq \mathbb{R}$ is a closed support. Let us assume that the first $2N+1$ integer moments of $f$, defined by 
\begin{equation}
M_k :=  \int_{\mathbb{S}} f(v) v^k \, \mathrm{d}v \quad \text{for $k \in \{0, 1, \dots, 2N \}$,}
\end{equation}
are finite and known. 
The objective of Multi-Gaussian KBMM  \cite{ctr10b}, which we refer here as EQMOM \cite{yuan2012} is to provide a continuous approximation $f^\app > 0$ defined such that
\begin{equation}\label{eq:quad}
M_k =  M^\app_k:= \int_{\mathbb{S}} f^\app(v) v^k \, \mathrm{d}v \quad \text{for $k \in \{0, 1, \dots, 2N \}$}.
\end{equation}
More precisely, $f^\app$ will be assumed to have the form:
\begin{equation}
f^\app(v) = \sum_{\alpha=1}^N \rho_\alpha \delta_\sigma(v,v_\alpha)
\end{equation}
where the given kernel density function (KDF) $\delta_\sigma(v,v_\alpha)$ ``tends'' to a Dirac delta function $\delta_{v_\alpha}(v)$ when $\sigma$ tends to zero and the $N$ weights $\rho_\alpha$, the $N$ abscissas  $v_\alpha$ and the ``spread" parameter $\sigma>0$ are determined from the first $2N+1$ integer moments of $f^\app$ by (\ref{eq:quad}).
Let us remark that, for a given type of KDF, it is not always possible to find such parameters. However, if the moment vector is not on the boundary of the moment space, one can always find a set of parameters such that (\ref{eq:quad}) holds for $k<2N$ and $M^\app_{2N}$ is as close as possible of $M_{2N}$.
In the following, we briefly describe Gaussian-EQMOM for $\mathbb{S}=\mathbb{R}$ and the associated algorithm for computing $f^\app$ from the moment set $\{M_0, \dots , M_{2N} \}$.

\subsection{Definition of Gaussian-EQMOM}\label{geqmom}

In Gaussian-EQMOM, $\delta_\sigma(v,v_\alpha)$ is chosen as a Gaussian PDF of variance $\sigma$ and centered at $v_\alpha$. 
The approximate NDF is then
\begin{equation}\label{Ngaus}
f^\G(v) := \sum_{\alpha = 1}^{N} \frac{\rho_{\alpha}}{\sigma \sqrt{2 \pi}} \exp \left( -\frac{\left( v-v_{\alpha} \right)^2}{2 \sigma^2} \right)
\end{equation}
and its moments are denoted $M^{\G}_k$.
In the case $\sigma = 0$,
$f$ is a weighted sum of $n \le N$ delta functions (i.e., an $n$-point distribution on a finite support). 
This case can be treated with an adaptive algorithm, described in \cite{yuan2011}, using the Chebyshev algorithm \cite{wheeler1974}.  
Thus, our focus in this work is on cases with $\sigma > 0$, which occur when $f$ is a continuous NDF.

Using the following definitions:
\begin{equation}\label{Nqmom}
\mu_k :=  \frac{1}{\sqrt{2 \pi}} \int_{\mathbb{R}} s^{k} e^{- s^2 / 2}  \, \mathrm{d}s ,  \quad M^*_k := \sum_{\alpha = 1}^{N} \rho_{\alpha} v_{\alpha}^k ,
\end{equation}
(note that $\mu_k = 0$ if $k$ is odd, otherwise $\mu_k \ge 1$), the nonnegative integer moments of $f^{G}$ can be written as
\begin{equation}\label{Nmom1}
M_k^{G} = \sum_{k_1 = 0}^{k} \binom{k}{k_1} \sigma^{k-k_1} \mu_{k-k_1} M^*_{k_1} ,
\end{equation}
which is a lower-triangular linear system of the form $M^G = A(\sigma) M^*$ where $A$ has unit eigenvalues.  Given $\sigma$ and the moment vector $M^G$, (\ref{Nmom1}) can thus be inverted to find the moment vector $M^*$.  This observation leads to the following moment-inversion algorithm.

\subsection{Moment-inversion algorithm for Gaussian-EQMOM}\label{moment-Gmom}

Given the moment set $\{M_0, \dots , M_{2N} \}$, the moment-inversion algorithm for Gaussian-EQMOM uses the following iterative procedure:  \\
Starting from $\sigma = 0$,
\begin{enumerate}
  \item compute $\{M^*_0, \dots , M^*_{2N-1} \}$ from (\ref{Nmom1}) by setting $M_k^G = M_k$.
  \item Use the Wheeler algorithm \cite{yuan2011} to compute $\{ \rho_1, \dots, \rho_N \}$ and $\{ v_1, \dots , v_N \}$.
  \item Compute $M^*_{2N}$ from (\ref{Nqmom}).
  \item Use (\ref{Nmom1}) to compute $M^G_{2N}$, compare with $M_{2N}$.
  \item If $M^G_{2N} < M_{2N}$, increase $\sigma$ and repeat from step 1 until $M^G_{2N} = M_{2N}$.
\end{enumerate}
In practice, this algorithm will converge if the moment set $\{M^*_0, \dots , M^*_{2N-1} \}$ is realizable for the given value of $\sigma$ used in step 1. Thus, if the moment set becomes unrealizable (which can be determined from the adaptive Wheeler algorithm \cite{yuan2011} in step 2), then the iterations are aborted and the largest value of $\sigma$ giving realizable moments is used (i.e., the value from the previous iteration).  By construction, this moment-inversion algorithm will yield $M_k^G = M_k$ for $k \in \{0, \dots, 2N-1 \}$ and $M_{2N}^G \le M_{2N}$ where the equality holds when the iterations are not aborted.  In \S\ref{s:quad-1D}, we consider the case with $N=2$ where an explicit formula is found for $\sigma$, and thus only steps 1 and 2 are required to compute the quadrature parameters.

\subsection{Evaluating integrals with Gaussian-EQMOM}\label{integral-1D}

In the kinetic equation describing multiphase flows, there are typically terms for collisions or breakage. In the moment transport equations, these terms lead to unclosed integrals of the form
\begin{equation} \label{eq:integral-1D}
\langle \mathcal{B} \rangle = \int_{\mathbb{R}}   \mathcal{B} (v) f(v) \, \mathrm{d}v , \quad
\langle \mathcal{C} \rangle = \int_{\mathbb{R}^2} \mathcal{C} (v_1,v_2) f(v_1) f(v_2) \, \mathrm{d}v_1 \mathrm{d}v_2
\end{equation}
where $\mathcal{B}$ and $\mathcal{C}$ are the breakage and collision kernels, respectively. Using the Gaussian-EQMOM representation for $f$, these integrals can be rewritten as
\begin{equation} \label{eq:integral-1Dmg}
\begin{gathered}
\langle \mathcal{B} \rangle = 
\sum_{\alpha=1}^{N} \rho_{\alpha} \frac{1}{\sqrt{\pi}}
\int_{\mathbb{R}} \mathcal{B} ( \sqrt{2} \sigma s + v_\alpha ) e^{-s^2} \, \mathrm{d}s , \\
\langle \mathcal{C} \rangle = 
\sum_{\alpha,\beta=1}^{N} \rho_{\alpha} \rho_{\beta} \frac{1}{\pi}
\int_{\mathbb{R}^2} \mathcal{C} ( \sqrt{2} \sigma s_1 + v_\alpha, \sqrt{2} \sigma s_2 + v_\beta ) e^{-(s_1^2+s_2^2)} \, \mathrm{d}s_1 \mathrm{d}s_2 .
\end{gathered}
\end{equation}
The remaining integrals in (\ref{eq:integral-1Dmg}) can then be approximated using an $M$-node Gauss-Hermite quadrature:
\begin{equation} \label{eq:integral-1Dmgh}
\begin{gathered}
\langle \mathcal{B} \rangle = 
\sum_{\alpha=1}^{N} \sum_{i=1}^{M} \rho_{\alpha} w_i \mathcal{B} ( \sqrt{2} \sigma s_i + v_\alpha ), \\
\langle \mathcal{C} \rangle = 
\sum_{\alpha, \beta =1}^{N} \sum_{i,j=1}^{M} \rho_{\alpha} \rho_{\beta} w_i w_j
\mathcal{C} ( \sqrt{2} \sigma s_i + v_\alpha, \sqrt{2} \sigma s_j + v_\beta )
\end{gathered}
\end{equation}
where $w_i$, $s_i$ are the corresponding Gauss-Hermite weights and abscissas \cite{gautschi04}. It is important to note that $M$ can be chosen independently from $N$, and thus that the accuracy of the quadrature approximation of (\ref{eq:integral-1Dmg}) by (\ref{eq:integral-1Dmgh}) does not depend on the number of moments.  The \emph{dual-quadrature form} of $f^G$ \cite{lage2011} is defined by
\begin{equation} \label{eq:dualquad}
f^G (v) =  \sum_{\alpha=1}^{N}  \sum_{i=1}^{M}  \rho_{\alpha i}  \delta ( v - \sqrt{2} \sigma s_i - v_\alpha ).
\end{equation}
where $\rho_{\alpha i} =  \rho_{\alpha} w_i $.

\vspace{12pt}

The reader can note that the Gaussian KDF used in (\ref{Ngaus}) could be replaced by any other normalized, symmetric function of $s = (v - v_{\alpha})/\sigma$. The moment-inversion algorithm in \S\ref{moment-Gmom} would remain unchanged, but it would be necessary to find the $M$-node Gaussian quadrature corresponding to the chosen KDF in order to efficiently evaluate the integrals in \S\ref{integral-1D}. Alternatively, \cite{ross2012} uses a two-node B-spline quadrature with a compact support $\mathbb{S} \subset \mathbb{R}$ for which the integral can be evaluated analytically. The B-spline quadrature can be extended to $N$ nodes using the methods described above.

\section{Application of 1-D two-node Gaussian-EQMOM to kinetic equations}\label{s:quad-1D}

We first introduce the two-node Gaussian-EQMOM for the NDF $f(t,x,v)$ in 1-D phase/real space for the kinetic equation:
\begin{equation} \label{equ:cinetique}
\partial_t f + v \partial_x f + \partial_v ( \mathcal{A} f ) = 0, \quad t> 0, \, x \in \mathbb{R}, \, v \in \mathbb{R}, 
\end{equation}
with initial condition $f(0,x,v) = f_0(x,v)$. The acceleration $\mathcal{A}$ is a real-valued function of $v$. The exact solution for free transport (when $\mathcal{A}=0$) is given by $f(t,x,v) = f(0,x-vt,v) = f_0(x-vt,v)$. In this work, we seek an approximation of $f(t,x,v)$ in the form of a two-node Gaussian-EQMOM with weights $\rho_1(t,x) > 0$, $\rho_2(t,x) > 0$, velocity abscissas $v_1(t,x)$, $v_2(t,x)$ and standard deviation $\sigma (t,x) \ge 0$. These five parameters are determined from the 1-D moment transport equations.

\subsection{1-D moment transport equations}\label{s:mom-1D}

Defining the $i^\text{th}$-order moment: 
\[
M_i (t,x) = \int_\mathbb{R} f(t,x,v) v^i \, \mathrm{d}v, \quad i=0, \dots ,K; \quad K \in \mathbb{N};
\]
the associated governing equations are easily obtained from (\ref{equ:cinetique}) after multiplication by $v^i$ and integration over $v$:
\[
\partial_t M_i + \partial_x M_{i+1} = \overline{\mathcal{A}}_i, \ i \geq 0,
\]
where the (unclosed)\footnote{The acceleration terms will be closed if $\mathcal{A}$ is affine: $\mathcal{A}(t,x,v) = - a(t,x) v + b(t,x)$, in which case the moment acceleration term can be written as $\overline{\mathcal{A}}_{k} = k (-a M_{k}+bM_{k-1})$. In gas-particle flows, this limit corresponds to Stokes drag in a stationary fluid.} moment acceleration term is
\begin{equation} \label{accel-1D}
\overline{\mathcal{A}}_i =  \int_\mathbb{R} i \mathcal{A}(v) f(t,x,v) v^{i-1} \, \mathrm{d}v .
\end{equation}
For simplicity, we will focus our attention on the five-moment model and its abstract form:
\begin{equation} \label{modele_bipic}
\begin{aligned}
\partial_t M_0 + \partial_x M_1 &= 0, \\
\partial_t M_1 + \partial_x M_2 &= \overline{\mathcal{A}}_1, \\
\partial_t M_2 + \partial_x M_3 &= \overline{\mathcal{A}}_2, \\
\partial_t M_3 + \partial_x M_4 &= \overline{\mathcal{A}}_3, \\
\partial_t M_4 + \partial_x \overline{M}_5 &= \overline{\mathcal{A}}_4.
\end{aligned}
\quad \Longrightarrow \quad
\partial_t {\bf M} + \partial_x {\bf F}({\bf M}) = \overline{\bfA},
\end{equation}
with $\bfM = (M_0,\dots,M_4)^t$, $\bfF ( \bfM ) = (M_1,\dots,M_4,\overline{M}_5)^t$ and $\overline{\bfA} = (0,\overline{\mathcal{A}}_1,\dots,\overline{\mathcal{A}}_4)^t$. This model is closed provided that $\overline{M}_5$ and $\overline{\bfA}$ are defined as functions of $\bfM$. Here we propose to define these functions using two-node Gaussian-EQMOM.

\subsection{Two-node Gaussian-EQMOM}\label{mia}

The function $f^\G$ has exact moments $M_i^\G$ of orders $i=0,...,5$ given by (\ref{Nmom1}).  The moment closure for (\ref{modele_bipic}) then naturally consists in setting $\overline{M}_5 =  M_5^\G$ where the five unknowns $\rho_1$, $\rho_2$, $v_1$, $v_2$ and $\sigma$ are found by solving the nonlinear system $M_i = M_i^\G, \ i=0,\dots, 4$; which is clearly equivalent to solving the system 
\begin{equation} \label{syst:quadrature_bis}
\begin{aligned}
M_0 &= \rho_1 + \rho_2, \\
M_1 &= \rho_1 v_1 + \rho_2 v_2, \\
M_2 - \sigma^2 M_0 &= \rho_1 v_1^2 + \rho_2 v_2^2, \\
M_3 - 3 \sigma^2 M_1 &= \rho_1 v_1^3 + \rho_2 v_2^3, \\
M_4 - 6 \sigma^2 M_2 + 3 \sigma^4 M_0 &= \rho_1 v_1^4 + \rho_2 v_2^4.
\end{aligned} 
\end{equation}
It remains to prove that this system is well-posed in the following proposition. 

\vspace{10pt}

\begin{proposition}[\textrm{Two-Node Gaussian-EQMOM}] \label{prop_quadrature}

For ${\bf M} = (M_0,M_1,M_2,M_3,M_4)^t$ such that $M_0 > 0$, define 
\[
e = \frac{M_0 M_2-M_1^2}{M_0^2}, \quad 
q = \frac{(M_3M_0^2-M_1^3)-3M_1(M_0M_2-M_1^2)}{M_0^3},
\]
and
\[
\eta = \frac{-3M_1^4+M_4M_0^3-4M_0^2M_1M_3+6M_0M_1^2M_2}{M_0^4}.
\]
System (\ref{syst:quadrature_bis}) is well-defined on the phase space $\Omega$ given by
\[
\Omega = 
\left\{ \bfM, \ M_0 > 0, \ e > 0, \ \eta > e^2 + \frac{q^2}{e}, \ \text{and} \  \eta \leq 3e^2 \ \text{if} \ q = 0 \right\}.
\]
Setting ${\bf U} = (\rho_1,\rho_2,\rho_1 v_1,\rho_2 v_2,\sigma)^t$, the function ${\bf U} = {\bf U}(\bfM)$ is one-to-one and onto when $v_1 \neq v_2$, and for all $v_1$ and $v_2$ provided that we set $\rho_1=\rho_2$ in the case $v_1 = v_2$. Moreover, $\sigma^2$ is given by the unique real root in $(0,e]$ of the third-order polynomial
\[
\begin{aligned}
\mathcal{P}(\sigma_0) &= 2 \sigma_0^3 + (\eta - 3 e^2) \sigma_0 + q^2, \\
\sigma_0 &= \sigma^2 - e. \\
\end{aligned}
\]
\end{proposition}

\emph{Proof}. 
By setting $\overline{\rho}_1 = \frac{\rho_1}{M_0}$, $\overline{\rho}_2 = \frac{\rho_2}{M_0}$, $\overline{v}_1 = v_1-\frac{M_1}{M_0}$, $\overline{v}_2 = v_2-\frac{M_1}{M_0}$, solving (\ref{syst:quadrature_bis}) is equivalent to solving
\[
\begin{aligned}
1 &= \overline{\rho}_1 + \overline{\rho}_2, \\
0 &= \overline{\rho}_1  \overline{v}_1 + \overline{\rho}_2  \overline{v}_2, \\
e - \sigma^2 &= \overline{\rho}_1 \overline{v}_1^2 + \overline{\rho}_2  \overline{v}_2^2, \\
q &= \overline{\rho}_1 \overline{v}_1^3 + \overline{\rho}_2 \overline{v}_2^3, \\
\eta - 6 \sigma^2 e + 3 \sigma^4 &= \overline{\rho}_1 \overline{v}_1^4 + \overline{\rho}_2  \overline{v}_2^4 ,
\end{aligned}
\]
with $e = (M_0 M_2-M_1^2)/M_0^2$, $q = ((M_3M_0^2-M_1^3)-3M_1(M_0M_2-M_1^2))/M_0^3$, $\eta = (-3M_1^4+M_4M_0^3-4M_0^2M_1M_3+6M_0M_1^2M_2)/M_0^4$. Dropping the overlines for the sake of clarity, it is then a matter of uniquely solving the following nonlinear system in $(\rho_1,\rho_2,v_1,v_2,\sigma^2)$:
\begin{equation} \label{sys_bigauss_centre}
\begin{aligned}
\rho_1 + \rho_2 &= 1, \\
\rho_1 v_1 + \rho_2 v_2 &= 0, \\
\rho_1 v_1^2 + \rho_2 v_2^2 &= e - \sigma^2, \\
\rho_1 v_1^3 + \rho_2 v_2^3 &= q, \\
\rho_1 v_1^4 + \rho_2 v_2^4 &= \eta - 6 \sigma^2 e + 3 \sigma^4.
\end{aligned}
\end{equation}
For any given value of $\sigma$ such that  $e > \sigma^2$ (or $e \geq \sigma^2$ if $q=0$), it is proved in \cite{fox08b} that the first four equations allow us to find $(\rho_1,\rho_2,v_1,v_2)$. We will then focus on the last equation to find $\sigma^2$.

In the case $q=0$, the second and fourth equations yield $\rho_1 v_1 (v_1^2 - v_2^2) = 0$, and $\rho_2 v_2 (v_1^2 - v_2^2) = 0$, which gives $v:=v_1=-v_2$. We then get $\rho_1=\rho_2=1/2$ and $\sigma^2 = e - v^2$, $2 v^4 = 3 e^2 - \eta$. Recall that our objective is now to uniquely determine $v$ and $\sigma > 0$ such that $\sigma^2 \leq e$. A necessary and sufficient condition is then clearly $\eta \in (e^2,3e^2)$. Note that the case $\eta=e^2$ would lead to $\sigma=0$, meaning that the Gaussian functions degenerate into two Dirac delta functions that correspond to the usual quadrature. The case $\eta=3e^2$ gives $v=0$ and both Gaussian functions coincide.

In the case $q\neq0$, from (\ref{sys_bigauss_centre}), we observe by using the usual algebra of quadrature methods and by setting $\sigma_0 = v_1v_2$ and $\sigma_1=-(v_1+v_2)$, that
\[
\begin{aligned}
e-\sigma^2+\sigma_0 &= 0, \\
q + \sigma_1 (e-\sigma^2) &= 0, \\
\eta - 6 \sigma^2 (e-\sigma^2) + 3 \sigma^4 + \sigma_1q + \sigma_0 (e-\sigma^2) &= 0.
\end{aligned}
\]
The last equation then gives that $\sigma_0=\sigma^2-e$ is a root of the third-order polynomial $\mathcal{P}(\sigma_0) = 2 \sigma_0^3 + (\eta - 3 e^2) \sigma_0 + q^2$. Note that one must have $\sigma_0\in (-e,0)$ to fulfill the condition $e> \sigma^2$ and to be able to reconstruct $\sigma > 0$ from $\sigma_0$. First, since $\lim_{\sigma_0 \to - \infty} \mathcal{P} = - \infty$, $\mathcal{P}(0) > 0$ and $\mathcal{P}^{''}(\sigma_0) = 12 \sigma_0$, there exists a unique root $\sigma_0 < 0$ of $\mathcal{P}$. It then follows that $\sigma_0 > -e$ if and only if $\mathcal{P}(-e) < 0$, that is if and only if $\eta > e^2 + q^2/e$. This concludes the proof.
\endproof

\subsection{Calculation of $\sigma^2$}

The three roots of $\mathcal{P}(\sigma_0)$ can be found analytically (one is real, and two are complex conjugates).  The real root yields 
\begin{equation} \label{eq:sigma}
\sigma^2 = e\left( 1 - c_3 + \frac{c_1}{c_3} \right)
\end{equation}
where
\[
c_1 = \frac{1}{6}\left( \frac{\eta}{e^2} - 3 \right),  \quad
c_2 = \frac{q^2}{4 e^{3}}   \quad \text{and} \quad
c_3 = \left[ \left( c_1^3 + c_2^2 \right)^{1/2} + c_2 \right]^{1/3} .
\]
Because $e$, $q$ and $\eta$ depend only on the moments and not on the weights and abscissas, $\sigma^2$ can be computed directly and thus the right-hand side of (\ref{sys_bigauss_centre}) is known.  

In transition zones between nondegenerate cases with $v_1 \neq v_2$ and the degenerate case with $v_1=v_2$, the two-node Gaussian-EQMOM based on (\ref{eq:sigma}) can generate quadrature points with high velocity but negligible density. As the CFL condition is based on the maximum velocity, this behavior will drastically reduce the time step.
 To avoid the generation of high velocities, a regularization procedure is employed. This kind of strategy is often used for moment systems close to singular distributions, see for instance \cite{junk1998,schaerer2015}.
  The difference between the mean velocity of each Gaussian is given by $\Delta u = |u_2-u_1| = q/(e-\sigma^2)$, and thus is an increasing function of $\sigma^2$. If $\sigma^2=0$, the velocity difference will be bounded \cite{kah11_cms}. Regularization consists of using a limiter on $\sigma^2$ to control the maximum velocity difference:
\begin{equation}\label{eq:sigmar}
\sigma^2 = 
\begin{cases}
\text{$\sigma^2$ from (\ref{eq:sigma})}           & \text{if $\Delta u < \Delta u_\text{lim}$} \\
e - \dfrac{|q|}{l} & \text{if $\Delta u \ge \Delta u_\text{lim}$}
\end{cases}
\end{equation}
where
\begin{equation}
l = \Delta u_\text{lim} + \left( \Delta u_\text{max}-\Delta u_\text{lim}\right) \tanh \left( \dfrac{\Delta u-\Delta u_\text{lim}}{\Delta u_\text{max}-\Delta u_\text{lim}} \right).
\end{equation}
The user-defined parameters $\Delta u_\text{lim}$ and $\Delta u_\text{max}>\Delta u_\text{lim}$ control the velocity difference starting from which we begin to limit $\sigma$ and the maximum velocity difference, respectively. In practice, these parameters should be close to the local fluid-phase velocity. The difference between $\Delta u_\text{lim}$ and $\Delta u_\text{max}$ permits a smooth transition between regularized and non-regularized zones. Note that in the regularized zones only the fourth-order moment is not satisfied (i.e., $M^G_4 < M_4$).

For $M_0 >0$, the two-node Gaussian-EQMOM algorithm consists of the following three steps:
\begin{remunerate}
\item Given moments $\mathbf{M}$ in $\Omega$, compute $e$, $q$ and $\eta$.
\item Compute $\sigma^2$ from (\ref{eq:sigmar}).
\item Solve (\ref{sys_bigauss_centre}) according to \cite{fox08b} to find $\rho_1$, $\rho_2$, $v_1$, $v_2$ (i.e., variables without the overlines). In the case where $\sigma^2 = e$, set $\rho_2=v_2=0$ and $\rho_1=M_0$, $v_1 = M_1/M_0$.
\end{remunerate}
For $M_0=0$, set $\rho_1=\rho_2=0$ and (without loss of generality) $\sigma = v_1=v_2=0$. 

\subsection{Effect of the $\sigma^2$ limiter on the flux}

\begin{figure}[htb]
\begin{center}
\includegraphics[width = 0.47\textwidth]{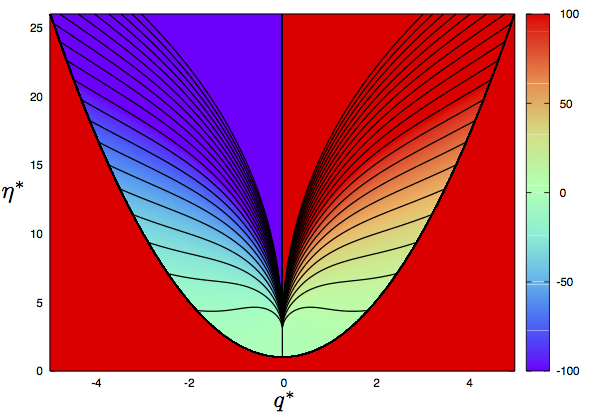}
\includegraphics[width = 0.47\textwidth]{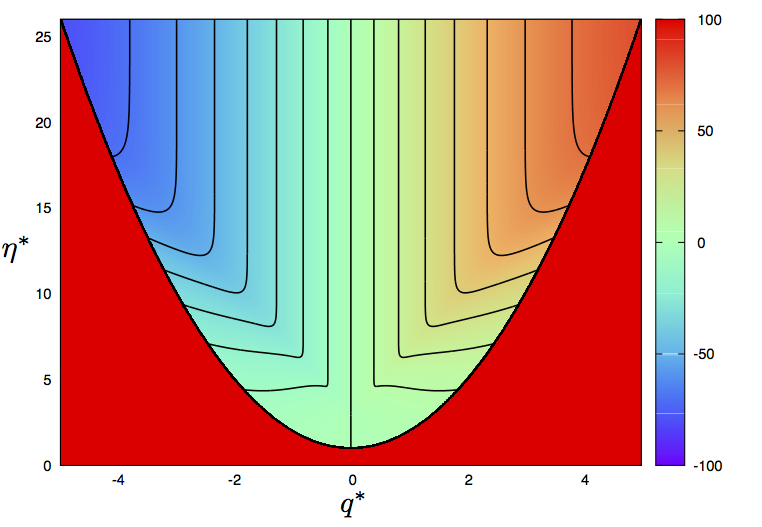}
\caption{Normalized flux $s^*$ without (left) and with (right) the limiter on $\sigma^2$.}
\end{center}
\label{fig:flux}
\end{figure}

When using a five-moment method, the moment space can be described through a set of normalized moments defined by Junk \cite{junk1998}.
Then, introducing the three central moments $e$, $q$, $\eta$ defined in Proposition~\ref{prop_quadrature} and a fourth one, occurring in the flux:
$$
s = \frac{M_5}{M_0} - 5 \frac{M_4M_1}{M_0^2} + 10 \frac{M_3M_1^2}{M_0^3}  -10 \frac{M_2M_1^3}{M_0^4} 
+4\frac{M_1^5}{M_0^5},
$$
the normalized moments, corresponding to the normalized NDF $f^*(w)=\frac{\sqrt e}{M_0}f\left(\sqrt e w +\frac{M_1}{M_0}\right)$,
are $(1,0,1,q^*,\eta^*,s^*)$ with $q^* = \frac{q}{e^{3/2}}$, $\eta^* = \frac{\eta}{e^2}$ and $s^*=\frac{s}{e^{5/2}}$.
In order to observe the effect of the $\sigma^2$ limiter on the flux, the normalized flux $s^*$ is plotted in Figure~\ref{fig:flux}, with and without the use of a limiter.
As observed for the closure with entropy maximization \cite{mcdonald2013}, without any limiter, the flux diverges quickly towards negative or positive infinity as the line $q^*=0$ with $\eta^*>3$ is approached.
However, with the limiter, the flux is smooth and does not diverge.

\subsection{Mathematical properties of moment system with two-node Gaussian-EQMOM}

The following theorem addresses an important mathematical property of system (\ref{modele_bipic}), which is its hyperbolicity.
This property was shown in the near thermodynamic equilibrium limit, characterized by $v_1\approx v_2$,  by Cheng and Rossmanith \cite{ross2012}. 
Here, the general case is addressed in the following theorem.
\begin{theorem}[\textrm{Hyperbolicity}] \label{th_hyp}
Assuming that the vector $\bfM=(M_0,M_1,M_2,M_3,M_4)^t$ lives in the space $\Omega$ defined in Proposition~\ref{prop_quadrature},  system (\ref{modele_bipic}) with the two-node Gaussian-EQMOM closure is hyperbolic.
\end{theorem}

\emph{Proof}. 
It is shown in Proposition~\ref{prop_quadrature} that one can define $U=(\rho_1,\rho_2,\rho_1v_1,\rho_2v_2,\sigma^2)^t$, the vector of the reconstruction variables, using the two-node Gaussian-EQMOM closure.  Moreover, $\sigma$ is not equal to zero since $\mathcal{P}(-e)=e\left(e^2+\frac{q^2}{e}-\eta\right)< 0$ with the notation of Proposition~\ref{prop_quadrature}.

The Jacobian matrix of system (\ref{modele_bipic}) is
$$
J = 
\begin{pmatrix}
0 & 1 & 0 & 0 & 0\\
0 & 0 & 1 & 0 & 0\\
0 & 0 & 0 & 1 & 0\\
0 & 0 & 0 & 0 & 1\\
\alpha & \beta & \gamma & \delta & \epsilon
\end{pmatrix}
$$
where  
$(\alpha,\beta,\gamma,\delta,\epsilon)=\dfrac{D\overline{M}_5}{DU}  . \left(\dfrac{D\bfM}{DU}\right)^{-1}$. After some tedious algebra, one obtains
$$
\begin{aligned}
\epsilon & = -2s_1+\bv,
\\
\delta & = -s_1^2-4s_0-2\bvd+10\sigma^2,
\\
\gamma & = -3s_0s_1-s_1\bvd+2s_0\bv+6\sigma^2(2s_1-\bv),
\\
\beta & = -2s_0\bvd -3s_0^2 +3\sigma^2(s_1^2+2\bvd+4s_0)-15\sigma^4,
\\
\alpha & = s_0^2\bv + \sigma^2(3s_0s_1+s_1\bvd-2s_0\bv)-3\sigma^4(2s_1-\bv)
\end{aligned}
$$
where $s_0$, $s_1$, $\bv$ and $\bvd$ are defined by
$$
s_0=v_1v_2, \quad 
s_1=-(v_1+v_2), \quad
\bv = \frac{\rho_1v_1+\rho_2v_2}{\rho_1+\rho_2} \quad \text{and} \quad
\bvd = \frac{\rho_1v_1^2+\rho_2v_2^2}{\rho_1+\rho_2} .
$$
The corresponding characteristic polynomial is 
$$
P(X)=-X^5+\epsilon\,X^4+\delta\,X^3+\gamma\,X^2+\beta\,X+\alpha.
$$
We then want to prove that this polynomial has five distinct roots.

In the case $v_1=v_2$ ($=\bv$), using the change of variables $X = Y +\bv$ leads to
$$
P(Y+\bv) = -Y^5+10\,\sigma^2 Y^3 -15\sigma^4 Y 
= -Y\,\left[Y^2-\left(5+\sqrt{10}\right)\sigma^2\right]\,\left[Y^2-\left(5-\sqrt{10}\right)\sigma^2\right] ,
$$
and hence there are five distinct roots when $\sigma \ne 0$.

In the case $v_1\ne v_2$, using the change of variables $X = \bv + \dfrac{v_2-v_1}{\rho_1+\rho_2}Z$ and defining $\Sigma= \dfrac{\rho_1+\rho_2}{v_2-v_1}\sigma$, leads to
$$
P\left(\bv+\frac{v_2-v_1}{\rho_1+\rho_2}Z\right) = -\left(\frac{v_2-v_1}{\rho_1+\rho_2}\right)^5 Q(Z)
$$
with
\begin{multline*}
Q(Z)=
Z^5 
+ 2(\rho_2-\rho_1) Z^4 
+ \left[ \rho_1^2-4\rho_1\rho_2+\rho_2^2-10\,\Sigma^2 \right] Z^3 
- 2(\rho_2-\rho_1)\left[ \rho_1\rho_2 +6 \Sigma^2 \right] Z^2  \\
+ \left[ \rho_1^2\rho_2^2-3(\rho_1^2-4\rho_1\rho_2+\rho_2^2)\Sigma^2 +15\Sigma^4 \right] Z
+ 2(\rho_2-\rho_1)\left[\rho_1\rho_2 +3\Sigma^2\right]\Sigma^2.
\end{multline*}
To prove that $Q(Z)$ has five distinct real roots, we use Sturm's theorem, which yields the number of distinct real roots in a given interval. 

Let us define the Sturm sequence of polynomials:
$P_0=Q$, $P_1 = Q' $ and, for any $n \in \{0,1,2,3\}$, $-P_{n+2}$ is the remainder of the Euclidean division of $P_{n+1}$ by $P_{n}$.
With the use of the formal calculation software Maxima, one can compute this sequence and remark that the coefficient of the highest-order term of each $P_n$ is positive for $n \in \{0,1,2,3,4,5\}$, since they are positive coefficient polynomial functions of $\rho_1$, $\rho_2$ and $\Sigma^2$.

Let $m(\xi)$ denote the number of sign changes (ignoring zeroes) in the sequence  $\{ P_0(\xi), P_1(\xi), \dots, P_5(\xi) \}$.  For $b$ large enough, one has $m(b)=0$ since each $P_n$ tends to $+\infty$ as $\xi \to +\infty$.  Likewise, for $a$ small enough, one has $m(a)=5$. Sturm's theorem then shows that $Q$ has five real distinct roots in $[a,b]$. This concludes the proof.
\endproof

\subsection{Kinetic-based flux}\label{flux}

In our numerical implementation to solve (\ref{modele_bipic}), the spatial fluxes ${\bf F}({\bf M})$ are computed using a kinetic-based definition:
\begin{equation}\label{fluxfunc}
F_i(t,x) = \int_{0}^{\infty} f(t,x,v) v^{i+1} \, \mathrm{d} v + \int_{-\infty}^0 f(t,x,v) v^{i+1} \, \mathrm{d} v, \quad i=0,\dots,4;
\end{equation}
where the decomposition into positive and negative directions is used to define the flux function as proposed in \cite{Bouchut03}. The numerical representation of the flux function is a critical point in moment transport methods \cite{vikas2011} because only \emph{realizable} moment sets can be successfully inverted.  Formally, we close (\ref{fluxfunc}) using the two-node Gaussian-EQMOM:
\begin{equation}\label{fluxfunc1a}
F_i(t,x) = \sum_{\alpha=1}^{2} \rho_{\alpha} \left[ \langle v^{i+1} \rangle^{+}_{\alpha} + \langle v^{i+1} \rangle^{-}_{\alpha}  \right], \quad i=0,\dots,4;
\end{equation}
where
\begin{equation}\label{halfmom}
\begin{gathered}
\langle v^{i} \rangle^{+}_{\alpha} := \frac{1}{\sqrt{\pi}} \int_{\frac{-v_\alpha}{\sqrt{2} \sigma}}^{\infty}  ( v_\alpha + \sqrt{2} \sigma s )^{i} e^{-s^2}\, \mathrm{d} s, \\
\langle v^{i} \rangle^{-}_{\alpha} := \frac{1}{\sqrt{\pi}} \int^{\infty}_{\frac{v_\alpha}{\sqrt{2} \sigma}}  ( v_\alpha - \sqrt{2} \sigma s )^{i} e^{-s^2} \, \mathrm{d} s,
\end{gathered}
\end{equation}
can be computed analytically. 
To design a first-order scheme, this decomposition is sufficient as it corresponds to an upwind scheme at the kinetic level.  
In this case, the transport scheme is shown to be realizable with the following CFL-like condition: $\frac{\Delta t}{\Delta x}\max_\alpha (|v_\alpha|+1.8 \sigma\sqrt 2)\le 1$ (see Appendix~\ref{s:realizability}).
For a quasi-high-order scheme \cite{vikas2011}, the spatial fluxes can be found from (\ref{fluxfunc1a}) by employing a high-order spatial reconstruction for $\rho_{\alpha}$ and a first-order reconstruction for the abscissas $v_{\alpha}$ and $\sigma$. 
In summary, the numerical fluxes are computed as follows:
\begin{remunerate}
\item Given moments $\mathbf{M}$, compute $\rho_1$, $\rho_2$, $v_1$, $v_2$, and $\sigma$ using the two-node Gaussian-EQMOM algorithm in \S\ref{mia}.
\item Compute kinetic-based fluxes from (\ref{fluxfunc1a}) using (\ref{halfmom}). 
\item Compute finite-volume numerical moment fluxes as described in \cite{vikas2011}.
\end{remunerate}
The case where $\sigma =0$ (i.e., $f$ is composed of two delta functions) is handled by the adaptive Wheeler algorithm \cite{yuan2011}.

\section{Extension of Gaussian-EQMOM to 2-D phase space}\label{s:gquad-2D}

Consider a 2-D phase space with NDF $f(\vv)$ for $\vv = (u, v)^t$ and define the bivariate moments
\begin{equation}
M_{i,j} := \int_{\mathbb{R}^2}  f(\vv) u^i v^j \, \mathrm{d} \vv, \ i,j =0, \dots ,K; \ K \in \mathbb{N} .
\end{equation}
Assuming that these moments exist, in the following we propose a bivariate extension of Gaussian-EQMOM using ideas from CQMOM \cite{yuan2011}.

\subsection{Definition of 2-D Gaussian-ECQMOM}\label{s:mgd-2D}

For clarity, we limit our discussion here to two-node quadrature in each direction of phase space (i.e.\ a total of four nodes).  Nevertheless, the same methodology can be used to develop the formulas for more nodes.  For the four-node quadrature, we define an approximate bivariate NDF by 
\begin{equation}\label{multigausa1}
f^{\G}_{12}(\vv) := 
\sum_{\alpha=1}^{2} \rho_{\alpha} g(u ; u_\alpha , \sigma_{1} ) 
\left( \sum_{\beta=1}^{2} \rho_{\alpha \beta} g(v - V(u) ; v_{\alpha \beta} , \sigma_{2 \alpha} ) \right)
\end{equation}
where the Gaussian KDF is  
\begin{equation}\label{gaussa1}
g( u ; \mu , \sigma ) := \frac{1}{\sigma \sqrt{2 \pi }} \exp \left( - \frac{( u - \mu )^2}{ 2 \sigma^2} \right) .
\end{equation}
The form in (\ref{multigausa1}) is based on $v$ conditioned on $u$.  An analogous form $f^{\G}_{21}(\vv)$ with $u$ conditioned on $v$ is found by permuting $u$ and $v$. For clarity, we describe the properties of Gaussian extended CQMOM (or Gaussian-ECQMOM) for $f^{\G}_{12}$.  

The function $V(u)$ in (\ref{multigausa1}) replaces the phase-space rotation introduced in \cite{fox08} and is defined to have the following properties:
\begin{gather}\label{vuprop1}
\sum_{\alpha=1}^{2} \rho_{\alpha} \int_{\mathbb{R}}  V(u) g(u ; u_\alpha , \sigma_{1} ) \sd u = M_{0,1} \\ \intertext{and} \label{vuprop2}
\sum_{\alpha=1}^{2} \rho_{\alpha} \int_{\mathbb{R}}  u V(u) g(u ; u_\alpha , \sigma_{1} ) \sd u = M_{1,1} .
\end{gather}
Here, we set\footnote{A more general choice that uses all of the moments is $V(u) = \sum_{n=0}^{4} a_n u^n$ with $a_n$ found from the linear system 
$$\sum_{n=0}^{4} \left[ \sum_{\alpha=1}^{2} \rho_{\alpha} \int_{\mathbb{R}}  u^{i+n} g(u ; u_\alpha , \sigma_{1} ) \sd u \right] a_n = M_{i,1} \quad \text{for $i \in \{ 0,1,2,3,4 \}$}.$$ 
This choice is valid if $\sigma_1 > 0$, but care must be taken if the coefficient matrix becomes poorly conditioned. Note that the quantity on the left-hand side in the square brackets is $M^G_{i+n,0}$, i.e., the univariate integer moments from two-node Gaussian-EQMOM. Thus the coefficient matrix is a Hankel matrix.}  $V(u) = a_0 + a_1 u$ where
\begin{align}
a_0 &= \frac{M_{2,0} M_{0,1} - M_{1,0} M_{1,1}}{M_{0,0} M_{2,0} - M_{1,0}^2} = \mu_v - \mu_u a_1 \\ \intertext{and}
a_1 &= \frac{M_{0,0} M_{1,1} - M_{1,0} M_{0,1}}{M_{0,0} M_{2,0} - M_{1,0}^2} = \frac{\rho \sigma_v}{\sigma_u}
\end{align}
with $\mu_u = M_{1,0}/M_{0,0}$, $\sigma_{u}^2 = (M_{2,0}/M_{0,0}) - \mu_u^2$, $\mu_v = M_{0,1}/M_{0,0}$, $\sigma_{v}^2 = (M_{0,2}/M_{0,0}) - \mu_v^2$ and
$\rho = [ (M_{1,1}/ M_{0,0}) - \mu_u \mu_v ]/ (\sigma_u \sigma_v )$ is the correlation coefficient. Note that $V(u)$ is well defined if the variance in the first direction $\sigma_u^2$ is nonzero. In fact $V(u)$ is the conditional expected value of $v$ given $u$ for a bivariate Gaussian distribution. To simplify the notation, we introduce the following definition:
\begin{equation}
\langle u^i V^j \rangle_{\alpha} := \int_{\mathbb{R}}  u^i V(u)^j g(u ; u_\alpha , \sigma_{1} ) \sd u .
\end{equation}
This integral involves Gaussian integer moments up to order $i+j$, which are known functions of $u_\alpha$ and $\sigma_{1}$.

The integer moments of (\ref{multigausa1}) are
\begin{equation}\label{multimoma1}
M^{\G}_{i,j} = 
\sum_{\alpha=1}^{2} \rho_{\alpha} \int_{\mathbb{R}}  u^i g(u ; u_\alpha , \sigma_{1} ) 
\left( \sum_{\beta=1}^{2} \rho_{\alpha \beta} \int_{\mathbb{R}}  v^j g(v - V(u) ; v_{\alpha \beta} , \sigma_{2 \alpha} ) \sd v  \right) \rd u
\end{equation}
or, by defining $v = y + V(u)$, as
\begin{equation}\label{multimoma2}
M^{\G}_{i,j} = 
\sum_{\alpha=1}^{2} \rho_{\alpha} \int_{\mathbb{R}} u^i g(u ; u_\alpha , \sigma_{1} ) 
\left( \sum_{\beta=1}^{2} \rho_{\alpha \beta} \int_{\mathbb{R}} [y + V(u) ]^j g(y ; v_{\alpha \beta} , \sigma_{2 \alpha} ) \sd y  \right) \rd u .
\end{equation}
For integer $j$, a binomial expansion leads to
\begin{equation}\label{multimoma3}
M^{\G}_{i,j} = 
\sum_{\alpha=1}^{2} \rho_{\alpha} \sum_{j_1=0}^{j} \binom{j}{j_1} \langle u^i V^{j-j_1} \rangle_\alpha \mu_{\alpha}^{j_1}
\end{equation}
where we have introduced the conditional moment of $(v-V)^j$ given $u = u_\alpha$ defined by
\begin{equation}\label{condmu}
\mu_{\alpha}^j := \sum_{\beta=1}^{2} \rho_{\alpha \beta} \int_{\mathbb{R}} y^j g(y ; v_{\alpha \beta} , \sigma_{2 \alpha} ) \sd y .
\end{equation}
It follows immediately from the properties of the Gaussian distribution and the definition of $\rho_{\alpha \beta}$ that $\mu_{\alpha}^0 = 1$ and $\mu_{\alpha}^1=0$.\footnote{$\mu_{\alpha}^0 = 1$ is the unique solution to (\ref{multimoma3}) with $i=0,1$ and $j=0$. Likewise, $\mu_{\alpha}^1=0$ is the solution for $i=0,1$ and $j=1$.  The latter makes use of the properties of $V(u)$ in (\ref{vuprop1}) and (\ref{vuprop2}).}  The form of (\ref{multimoma3}) leads to the following moment-inversion algorithm.

\subsection{Moment-inversion algorithm for 2-D Gaussian-ECQMOM}\label{mia-2D}

The first step uses the univariate moments $M_{i,0}$ with the algorithm in \S\ref{moment-Gmom} for 1-D Gaussian-EQMOM to find $\{ \rho_{1} , \rho_{2} ,  u_{1}, u_{2}, \sigma_{1}  \}$. There are two possible cases: (1) a nondegenerate case with $u_1 \neq u_2$, (2) a degenerate case with $u_1 = u_2$. The degenerate case occurs when the univariate moments $M_{i,0}$ are Gaussian.

\subsubsection{Nondegenerate case}

In order to determine the parameters $v_{\alpha \beta}$ and $\sigma_{2 \alpha}$, we must solve (\ref{multimoma3}) to find $\mu_{\alpha}^j$ for $j=2,3,4$. This requires six equations, which we define by taking $i=0,1$ and $j=2,3,4$.  It is straightforward to show that the conditional moments can be found from three $2 \times 2$ linear systems, which can be solved sequentially:
\begin{equation}\label{linsysa1}
\begin{aligned}
\sum_{\alpha=1}^{2} \rho_{\alpha} \mu_{\alpha}^2            &= M_{0,2} - \sum_{\alpha=1}^{2} \rho_{\alpha} \langle   V^{2} \rangle_\alpha  , \\
\sum_{\alpha=1}^{2} \rho_{\alpha} u_{\alpha} \mu_{\alpha}^2 &= M_{1,2} - \sum_{\alpha=1}^{2} \rho_{\alpha} \langle u V^{2} \rangle_\alpha ;
\end{aligned}
\end{equation}
\begin{equation}\label{linsysa2}
\begin{aligned}
\sum_{\alpha=1}^{2} \rho_{\alpha} \mu_{\alpha}^3            &= M_{0,3} - 3 \sum_{\alpha=1}^{2} \rho_{\alpha} \langle V \rangle_\alpha \mu_{\alpha}^2 
- \sum_{\alpha=1}^{2} \rho_{\alpha} \langle   V^{3} \rangle_\alpha   , \\
\sum_{\alpha=1}^{2} \rho_{\alpha} u_{\alpha} \mu_{\alpha}^3 &= M_{1,3} - 3 \sum_{\alpha=1}^{2} \rho_{\alpha} \langle u V \rangle_\alpha \mu_{\alpha}^2
- \sum_{\alpha=1}^{2} \rho_{\alpha} \langle u V^{3} \rangle_\alpha ;
\end{aligned}
\end{equation}
\begin{equation}\label{linsysa3}
\begin{aligned}
\sum_{\alpha=1}^{2} \rho_{\alpha} \mu_{\alpha}^4            &= M_{0,4} - 4 \sum_{\alpha=1}^{2} \rho_{\alpha} \langle V \rangle_\alpha \mu_{\alpha}^3
- 6 \sum_{\alpha=1}^{2} \rho_{\alpha} \langle V^2 \rangle_\alpha \mu_{\alpha}^2 - \sum_{\alpha=1}^{2} \rho_{\alpha} \langle   V^{4} \rangle_\alpha  , \\
\sum_{\alpha=1}^{2} \rho_{\alpha} u_{\alpha} \mu_{\alpha}^4 &= M_{1,4} - 4 \sum_{\alpha=1}^{2} \rho_{\alpha} \langle u V \rangle_\alpha \mu_{\alpha}^3
- 6 \sum_{\alpha=1}^{2} \rho_{\alpha} \langle u V^2 \rangle_\alpha \mu_{\alpha}^2 - \sum_{\alpha=1}^{2} \rho_{\alpha} \langle u V^{4} \rangle_\alpha .
\end{aligned}
\end{equation}
The final step is to use the algorithm in \S\ref{moment-Gmom} for 1-D Gaussian-EQMOM with the moment set $\{ 1 , 0 , \mu_{\alpha}^2 , \mu_{\alpha}^3 , \mu_{\alpha}^4 \}$  to find the quadrature parameters $\{ \rho_{\alpha 1} , \rho_{\alpha 2} ,  v_{\alpha 1}, v_{\alpha 2}, \sigma_{2 \alpha}  \}$ for $\alpha \in \{ 1,2 \}$.

\subsubsection{Degenerate case}\label{sec:degcase}

For the degenerate case, we define a 2-node Gaussian-ECQMOM by 
\begin{equation}\label{multigausa1d}
f^{\G}_{12}(\vv) = M_{0,0} g(u ; \mu_u , \sigma_{u} ) \left( \sum_{\alpha =1}^{2} \rho_{\alpha} g(v - V(u) ; v_{\alpha} , \sigma_{2} ) \right)
\end{equation}
that can further degenerate\footnote{If the bivariate moment set is (anisotropic) Gaussian, then the quadrature will degenerate to the exact distribution with $w_1=1$, $v_1=0$ and $\sigma_{2}^2 = (1 - \rho^2 ) \sigma_{v}^2$.} to a bivariate Gaussian with a full covariance matrix when $\rho_2=0$.  The moments for this case are
\begin{equation}\label{multimoma3d}
M^{\G}_{i,j} = M_{0,0} \sum_{j_1=0}^{j} \binom{j}{j_1} \langle u^i V^{j-j_1} \rangle \mu^{j_1}
\end{equation}
where
\begin{equation}
\langle u^i V^j \rangle := \int_{\mathbb{R}} u^i V(u)^j g(u ; \mu_u , \sigma_{u} ) \sd u
\end{equation}
with $\langle V \rangle = M_{0,1} / M_{0,0} = \mu_v$ and $\langle u V \rangle = M_{1,1} / M_{0,0}$. The conditional moments are defined by
\begin{equation}\label{condmud}
\mu^j := \sum_{\alpha=1}^{2} \rho_{\alpha} \int_{\mathbb{R}} y^j g(y ; v_{\alpha} , \sigma_{2} ) \sd y .
\end{equation}
By definition, $\mu^0 = 1$ and $\mu^1 = 0$. To complete the quadrature, we find $\mu^2$, $\mu^3$ and $\mu^4$ from (\ref{multimoma3d}) using $i=0$ and $j=2,3,4$:
\begin{gather}
\mu^2 = \frac{M_{0,2}}{M_{0,0}} - \langle V^2 \rangle \\
\mu^3 = \frac{M_{0,3}}{M_{0,0}} - \langle V^3 \rangle - 3 \langle V \rangle \mu^2 \\
\mu^4 = \frac{M_{0,4}}{M_{0,0}} - \langle V^4 \rangle - 6 \langle V^2 \rangle \mu^2 - 4 \langle V \rangle \mu^3
\end{gather}
where the moments $\langle V^j \rangle$ are known. The final step is to use the algorithm in \S\ref{moment-Gmom} for 1-D Gaussian-EQMOM with the moment set $\{ 1 , 0 , \mu^2 , \mu^3 , \mu^4 \}$  to find the quadrature parameters $\{ \rho_{1} , \rho_{2} ,  v_{1}, v_{2}, \sigma_{2}  \}$.

For the nondegenerate case, the moment-inversion algorithm described above is able to recover 13 of the 16 \emph{extended optimal moments}\footnote{If we define $V(u)= \sum_{n=0}^{4} a_n u^n$, then all 16 moments can be recovered.} \cite{fox09} defined by
\[
\begin{bmatrix}
M_{0,0} & M_{0,1} & M_{0,2}      & M_{0,3}       & M_{0,4} \\
M_{1,0} & M_{1,1} & M_{1,2}      & M_{1,3}       & M_{1,4} \\
M_{2,0} & M_{2,1} &              &               & \\
M_{3,0} & M_{3,1} &              &               & \\
M_{4,0} & M_{4,1} &              &               & \\
\end{bmatrix}   .
\]
The reader can note that the formulas developed in this section for four-node Gaussian-ECQMOM can be extended to $N \times N$ nodes in a relatively straightforward manner.  

\subsection{Evaluating integrals with 2-D Gaussian-ECQMOM}\label{integral-2D}

In 2-D phase space, the evaluation of integrals using Gaussian-ECQMOM is very similar to \S\ref{integral-1D}. Consider again the unclosed integral
\begin{equation} \label{eq:integral-2D}
\langle \mathcal{B} \rangle = \int_{\mathbb{R}^2} \mathcal{B} (\vv) f(\vv) \, \mathrm{d} \vv .
\end{equation}
Using the 2-D Gaussian-ECQMOM representation $f^G_{12}$, this integral can be rewritten as
\begin{equation} \label{eq:integral-2Dmg}
\langle \mathcal{B} \rangle = 
\sum_{\alpha , \beta=1}^{2} 
\frac{\rho_{\alpha} \rho_{\alpha \beta} }{\pi}
\int_{\mathbb{R}^2} 
\mathcal{B} ( \sqrt{2} \sigma_1 s_1 + u_{\alpha} , \sqrt{2} \sigma_{2 \alpha} s_2 + V( \sqrt{2} \sigma_1 s_1 + u_{\alpha}) + v_{\alpha \beta}  ) 
e^{- \ss^2}  \, \mathrm{d} \ss 
\end{equation}
where $\ss = (s_1 , s_2 )^t$. The remaining integrals in (\ref{eq:integral-2Dmg}) can then be approximated using an $M$-node Gauss-Hermite quadrature:
\begin{equation} \label{eq:integral-2Dmgh}
\langle \mathcal{B} \rangle = 
\sum_{\alpha, \beta =1}^{2} \rho_{\alpha} \rho_{\alpha \beta} \sum_{i,j=1}^{M} 
w_{i} w_{j} 
\mathcal{B} ( \sqrt{2} \sigma_1 s_{i} + u_{\alpha} , \sqrt{2} \sigma_{2 \alpha} s_{j} + V( \sqrt{2} \sigma_1 s_{i} + u_{\alpha} ) + v_{\alpha \beta}  )
\end{equation}
where $w_{i}$ and $s_{i}$ are the Gauss-Hermite weights and abscissas \cite{gautschi04}. The extension of (\ref{eq:integral-2Dmg}) to a 3-D (or higher) phase space is analogous to (\ref{eq:integral-2Dmgh}).  The \emph{dual-quadrature representation} is thus
\begin{equation}\label{eq:dualquad2d}
f^G_{12} (\vv) = 
\sum_{\alpha, \beta =1}^{2} \sum_{i,j=1}^{M} 
 \rho_{\alpha \beta i j} \delta ( u -  \sqrt{2} \sigma_1 s_{i} -u_{\alpha} ) \delta ( v - \sqrt{2} \sigma_{2 \alpha} s_{j} - V( u) - v_{\alpha \beta}  )  
\end{equation}
where $\rho_{\alpha \beta i j }  = \rho_{\alpha} \rho_{\alpha \beta}  w_{i} w_{j} $.

A second, and equally valid, evaluation of the integrals in (\ref{eq:integral-2D}) can be found using $f^G_{21}$. In general, if there exist $P$ permutations of the CQMOM conditioning variables (e.g., $P=2$ in 2-D and $P=6$ in 3-D), then each of them can be used to evaluate the integrals. The arithmetic average value of the $P$ permutations would represent the best overall estimate \cite{yuan2011}.

\section{Application of four-node Gaussian-ECQMOM to kinetic equations}\label{s:quad-2D}

Consider a 2-D velocity phase space with NDF $f(t,\xx,\vv)$ for $\xx = (x, y)^t$ and $\vv = (u, v)^t$ that satisfies the kinetic equation
\begin{equation} \label{equ:cinetique2}
\partial_t f + \vv \cdot \partial_\xx f + \partial_\vv \cdot ( \bfA f ) = 0, 
\quad t> 0, \, \xx \in \mathbb{R}^2, \, \vv \in \mathbb{R}^2, 
\end{equation}
with initial condition $f(0,\xx,\vv) = f_0(\xx,\vv)$. The acceleration $\bfA = (\mathcal{A}_x, \mathcal{A}_y)^t$ is a real-valued function of $\vv$. With a 2-D velocity phase space, we approximate the solution to $f$ using a Gaussian-ECQMOM for the bivariate moments. In this work, we will consider only the \emph{minimal} Gaussian-ECQMOM that uses four nodes in the 2-D velocity phase space. Nonetheless, the extension to more than four nodes would be analogous to the algorithm presented here. 

\subsection{2-D moment transport equations}\label{s:mom-2D}

Defining the bivariate moments
\[
M_{i,j}(t,\xx) = \int_{\mathbb{R}^2} f(t,\xx,\vv) u^i v^j \, \mathrm{d} \vv, \ i,j =0, \dots ,K; \ K \in \mathbb{N};
\]
the associated governing equations are easily obtained from (\ref{equ:cinetique2}):
\[
\partial_t M_{i,j} + \partial_x M_{i+1,j} + \partial_y M_{i,j+1} = \overline{\mathcal{A}}_{i,j}, \quad i,j \geq 0;
\]
where the (unclosed)\footnote{The acceleration terms will be closed if $\bfA$ is a linear function of the form $(a u, a v)^t$, in which case the moment acceleration term can be written as $\overline{\mathcal{A}}_{i,j} = - a (i + j) M_{i,j}$. In gas-particle flows, this limit corresponds to Stokes drag in a stationary fluid.} moment acceleration term is defined by
\begin{equation} \label{modele_accel-2D}
\overline{\mathcal{A}}_{i,j} = 
- \int_{\mathbb{R}^2} i \mathcal{A}_x(\vv) f(t,\xx,\vv) u^{i-1} v^{j} \, \mathrm{d} \vv 
- \int_{\mathbb{R}^2} j \mathcal{A}_y(\vv) f(t,\xx,\vv) u^{i} v^{j-1} \, \mathrm{d} \vv .
\end{equation}
We will consider in this work a 16-moment model:
\begin{equation} \label{modele_bipic2}
\begin{gathered}
\begin{aligned}
\partial_t M_{0,0} + \partial_x M_{1,0} + \partial_y M_{0,1} &= 0, \\
\partial_t M_{1,0} + \partial_x M_{2,0} + \partial_y M_{1,1} &= \overline{\mathcal{A}}_{1,0}, \\
\partial_t M_{0,1} + \partial_x M_{1,1} + \partial_y M_{0,2} &= \overline{\mathcal{A}}_{0,1}, \\
\partial_t M_{2,0} + \partial_x M_{3,0} + \partial_y M_{2,1} &= \overline{\mathcal{A}}_{2,0}, \\
\partial_t M_{1,1} + \partial_x M_{2,1} + \partial_y M_{1,2} &= \overline{\mathcal{A}}_{1,1}, \\
\partial_t M_{0,2} + \partial_x M_{1,2} + \partial_y M_{0,3} &= \overline{\mathcal{A}}_{0,2}, \\
\partial_t M_{3,0} + \partial_x M_{4,0} + \partial_y M_{3,1} &= \overline{\mathcal{A}}_{3,0}, \\
\partial_t M_{2,1} + \partial_x M_{3,1} + \partial_y \overline{M}_{2,2} &= \overline{\mathcal{A}}_{2,1},
\end{aligned}
\quad
\begin{aligned}
\partial_t M_{1,2} + \partial_x \overline{M}_{2,2} + \partial_y M_{1,3} &= \overline{\mathcal{A}}_{1,2}, \\
\partial_t M_{0,3} + \partial_x M_{1,3}            + \partial_y M_{0,4} &= \overline{\mathcal{A}}_{0,3}, \\
\partial_t M_{4,0} + \partial_x \overline{M}_{5,0} + \partial_y M_{4,1} &= \overline{\mathcal{A}}_{4,0}, \\
\partial_t M_{3,1} + \partial_x M_{4,1} + \partial_y \overline{M}_{3,2} &= \overline{\mathcal{A}}_{3,1}, \\
\partial_t M_{1,3} + \partial_x \overline{M}_{2,3} + \partial_y M_{1,4} &= \overline{\mathcal{A}}_{1,3}, \\
\partial_t M_{0,4} + \partial_x M_{1,4} + \partial_y \overline{M}_{0,5} &= \overline{\mathcal{A}}_{0,4}, \\
\partial_t M_{4,1} + \partial_x \overline{M}_{5,1} + \partial_y \overline{M}_{4,2} &= \overline{\mathcal{A}}_{4,1}, \\
\partial_t M_{1,4} + \partial_x \overline{M}_{2,4} + \partial_y \overline{M}_{1,5} &= \overline{\mathcal{A}}_{1,4},
\end{aligned}
\end{gathered}
\end{equation}
which requires a closure for the fourth-order moment $\overline{M}_{2,2}$, the four fifth-order moments $\overline{M}_{5,0}, \overline{M}_{3,2}, \overline{M}_{2,3}, \overline{M}_{0,5}$, the four sixth-order moments $\overline{M}_{5,1}, \overline{M}_{4,2}, \overline{M}_{2,4}, \overline{M}_{1,5}$, and the acceleration terms. We propose to define these closures by reconstructing $f$ with four-node Gaussian-ECQMOM. If unclosed, the acceleration term $\overline{\bfA}$ can be evaluated using the dual-quadrature forms of $f^G_{12}$ and $f^G_{21}$, and taking the arithmetic average.  In our numerical examples, Stokes drag is used so that $\overline{\bfA}$ is closed in terms of the transported moments, and operator splitting is used for the fluxes and the acceleration.

\subsection{Mathematical properties of 2-D moment system with four-node Gaussian-ECQMOM }

Because we use a dimensional splitting to solve (\ref{modele_bipic2}), let us consider the transport part of the system in the $x$-direction for the moments used in the Gaussian-ECQMOM reconstruction conditioned on the $u$ velocity:
\begin{equation}\label{eq:systx}
\partial_t \bfM + \partial_x  \bfF( \bfM ) = 0
\end{equation}
with
$$
\bfM=(M_{0,0},M_{1,0},M_{2,0},M_{3,0},M_{4,0},
M_{0,1},M_{1,1},
M_{0,2},M_{1,2},
M_{0,3},M_{1,3},
M_{0,4},M_{1,4})^t
$$
and
$$
\bfF( \bfM )=(M_{1,0},M_{2,0},M_{3,0},M_{4,0},\overline{M}_{5,0},
M_{1,1},\overline{M}_{2,1},
M_{1,2},\overline{M}_{2,2},
M_{1,3},\overline{M}_{2,3},
M_{1,4},\overline{M}_{2,4})^t .
$$
Let us remark that the moments $M_{2,1}$, $M_{3,1}$ and $M_{4,1}$ are not considered here since their equations are redundant with those of (\ref{eq:systx}). 
Moreover, for the sake of simplicity, we use the reconstruction with $V(u)=0$ in (\ref{multigausa1}) to prove the hyperbolicity of system (\ref{eq:systx}).

\begin{theorem}[\textrm{Hyperbolicity}] \label{th_hyp2D}
Assuming that the moment-inversion algorithm for 2-D Gaussian-ECQMOM with $V(u)=0$ for the vector $\bfM$ is nondegenerate, system (\ref{eq:systx}) with this closure is hyperbolic.
\end{theorem}

\emph{Proof}. 
The Jacobian matrix $J_{2D}$ of the flux is block triangular:
$$
J_{2D} = 
\begin{pmatrix}
J     & 0 & 0 & 0 & 0\\
X_1 & A & 0 & 0 & 0\\
X_2 & 0 & A & 0 & 0\\
X_3 & 0 & 0 & A & 0\\
X_4 & 0 & 0 & 0 & A
\end{pmatrix}
$$ 
where $J$ is the $5\times 5$ Jacobian matrix given in Theorem \ref{th_hyp}, corresponding to the 1-D system. The $X_i$ are $2\times 5$ matrices and $A$ is a  $2\times 2$ matrix given by
$$
A=
\begin{pmatrix}
0 & 1\\
\nu & \xi
\end{pmatrix}
$$
with
\[
(\nu,\xi)=\dfrac{D\overline{M}_{2,1}}{D(M_{0,1},M_{1,1})}=\dfrac{D\overline{M}_{2,2}}{D(M_{0,2},M_{1,2})}=
\dfrac{D\overline{M}_{2,3}}{D(M_{0,3},M_{1,3})}=\dfrac{D\overline{M}_{2,4}}{D(M_{0,4},M_{1,4})}
\]
given by $\nu = -u_1u_2 + \sigma_1^2$ and $\xi = u_1+u_2$.

From Theorem \ref{th_hyp}, the matrix $J$ is diagonalizable with five distinct eigenvalues. Matrix $A$ is also diagonalizable with the following two eigenvalues:
$$
\lambda_{\pm} = \frac{u_1+u_2}{2} \pm \sqrt{\frac{(u_2-u_1)^2}{4}+\sigma_1^2} .
$$
Moreover, neither of these eigenvalues is an eigenvalue of $J$. It follows from the block-diagonal structure of the submatrix, found by eliminating the first five rows and columns of $J_{2D}$, that $J_{2D}$ is diagonalizable. This concludes the proof.
\endproof

\subsection{Kinetic-based flux}\label{flux2}

The moment transport system (\ref{modele_bipic2}) has the form
\[
\partial_t \bfM + \partial_\xx \cdot \bfF( \bfM ) = \overline{\bfA}
\]
with flux vector $\bfF =( \bfF_{x}, \bfF_{y} )^t$ for the 16-moment vector $\bfM$. In our numerical implementation, the components of the fluxes for moment $M_{i,j}$ are computed using a kinetic-based definition:
\begin{equation}\label{flux-1Dx}
F_{x;i,j} = \int_\mathbb{R} \left( \int_{0}^{\infty} f(t,\xx,\vv) u^{i+1} v^j \, \mathrm{d}u \right) \mathrm{d}v 
          + \int_\mathbb{R} \left( \int_{-\infty}^0 f(t,\xx,\vv) u^{i+1} v^j \, \mathrm{d}u \right) \mathrm{d}v ,
\end{equation}
\begin{equation}\label{flux-1Dy}
F_{y;i,j} = \int_\mathbb{R} \left( \int_{0}^{\infty} f(t,\xx,\vv) u^{i} v^{j+1} \, \mathrm{d}v \right) \mathrm{d}u 
          + \int_\mathbb{R} \left( \int_{-\infty}^0 f(t,\xx,\vv) u^{i} v^{j+1} \, \mathrm{d}v \right) \mathrm{d}u.
\end{equation}
Here, we describe how $F_{x;i,j}$ is computed by conditioning on the $u$ component of velocity.  The treatment of $F_{y;i.j}$ is done analogously by conditioning on the $v$ component of velocity.

Assuming that 1-D Gaussian-EQMOM has been applied to the $u$ direction, we can assume that $\rho_{\alpha}$, $u_{\alpha}$ and $\sigma_1$ are known. Thus, for the nondegenerate case with $u_1 \neq u_2$, (\ref{flux-1Dx}) can be written as
\begin{equation}\label{fluxfunc1a2D}
F_{x;i,j} = \sum_{\alpha=1}^{2} \rho_{\alpha} \sum_{j_1=0}^{j} \binom{j}{j_1}
\left[ \langle u^{i+1} V^{j-j_1} \rangle^{+}_{\alpha} + \langle u^{i+1} V^{j-j_1} \rangle^{-}_{\alpha}  \right] \mu_\alpha^{j_1}
\end{equation}
where $\mu_\alpha^j$ are the conditional moments of $(v-V)^j$ given $u_\alpha$, and
\begin{equation}\label{halfmom2D}
\begin{gathered}
\langle u^{i} V^j \rangle^{+}_{\alpha} := 
\frac{1}{\sqrt{\pi}} \int_{\frac{-u_\alpha}{\sqrt{2} \sigma_1}}^{\infty}  
( u_\alpha + \sqrt{2} \sigma_1 s )^{i} 
V( u_\alpha + \sqrt{2} \sigma_1 s )^{j}
e^{-s^2} \, \mathrm{d} s, \\
\langle u^{i} V^j \rangle^{-}_{\alpha} := 
\frac{1}{\sqrt{\pi}} \int^{\infty}_{\frac{u_\alpha}{\sqrt{2} \sigma_1}}  
( u_\alpha - \sqrt{2} \sigma_1 s )^{i}
V( u_\alpha - \sqrt{2} \sigma_1 s )^{j}
e^{-s^2} \, \mathrm{d} s ,
\end{gathered}
\end{equation}
can be computed analytically.  For the degenerate case ($\rho_2=0$), the same expressions are used but with the conditional moments $\mu^j$ defined in \S\ref{sec:degcase}.

\section{Numerical examples}\label{s:examples}

As example applications, we consider a Riemann problem with 1-D velocity phase space and a frozen turbulence problem with a 2-D velocity phase space. In each case, we solve the moment transport equations in (\ref{modele_bipic}) and (\ref{modele_bipic2}), respectively.

\subsection{1-D Riemann problem}\label{s:num-1D}

The initial conditions are defined on the real line with a step in the mean $u$ velocity at $x=0$: 
$$
U_\text{m} = \frac{M_1}{M_0} = 
\begin{cases} 
1 & \text{if $x < 0$,} \\
-1 &  \text{otherwise.}
\end{cases} 
$$
For all $x$, the initial density is unity and the velocity distribution function is Maxwellian with energy $\sigma^2 = 1/3$. The velocity distribution is assumed initially to be in equilibrium (i.e., $e_0=\sigma^2$). However, the discontinuous nature of the mean particle velocity quickly leads to particle trajectory crossing and a strongly non-equilibrium velocity distribution function.  

For 1-D phase space, a measure of the degree of non-equilibrium is the ratio $\sigma^2/e$, which is unity for an equilibrium distribution and zero when the distribution is composed entirely of Dirac delta functions.  In order to identify clearly deviations of the higher-order moments from their equilibrium values, we will use the following normalized moments: $e^* = e/e_0$ (energy), $q^* = q/e^{3/2}$ (skewness), $\eta^* = \eta/e^2$ (kurtosis); whose equilibrium values are $e^*=1$, $q^*=0$, and $\eta^*=3$.

In order to solve the moment equations numerically, the 1-D computational domain $-2 < x < 2$ is discretized into 402 finite-volume cells.  The spatial fluxes are treated using the first-order kinetic-based approach. The time step is chosen based on the largest magnitude of the abscissas $v_\alpha$ used to define the spatial fluxes with a CFL number of 0.5. For comparison, results are shown for $\sigma^2$ found with (\ref{eq:sigma}) and with the velocity limiter in (\ref{eq:sigmar}).


\begin{figure}[htb]
\begin{center}
\includegraphics[width = 1.0\textwidth,clip=true, trim= 40 0 40 0]{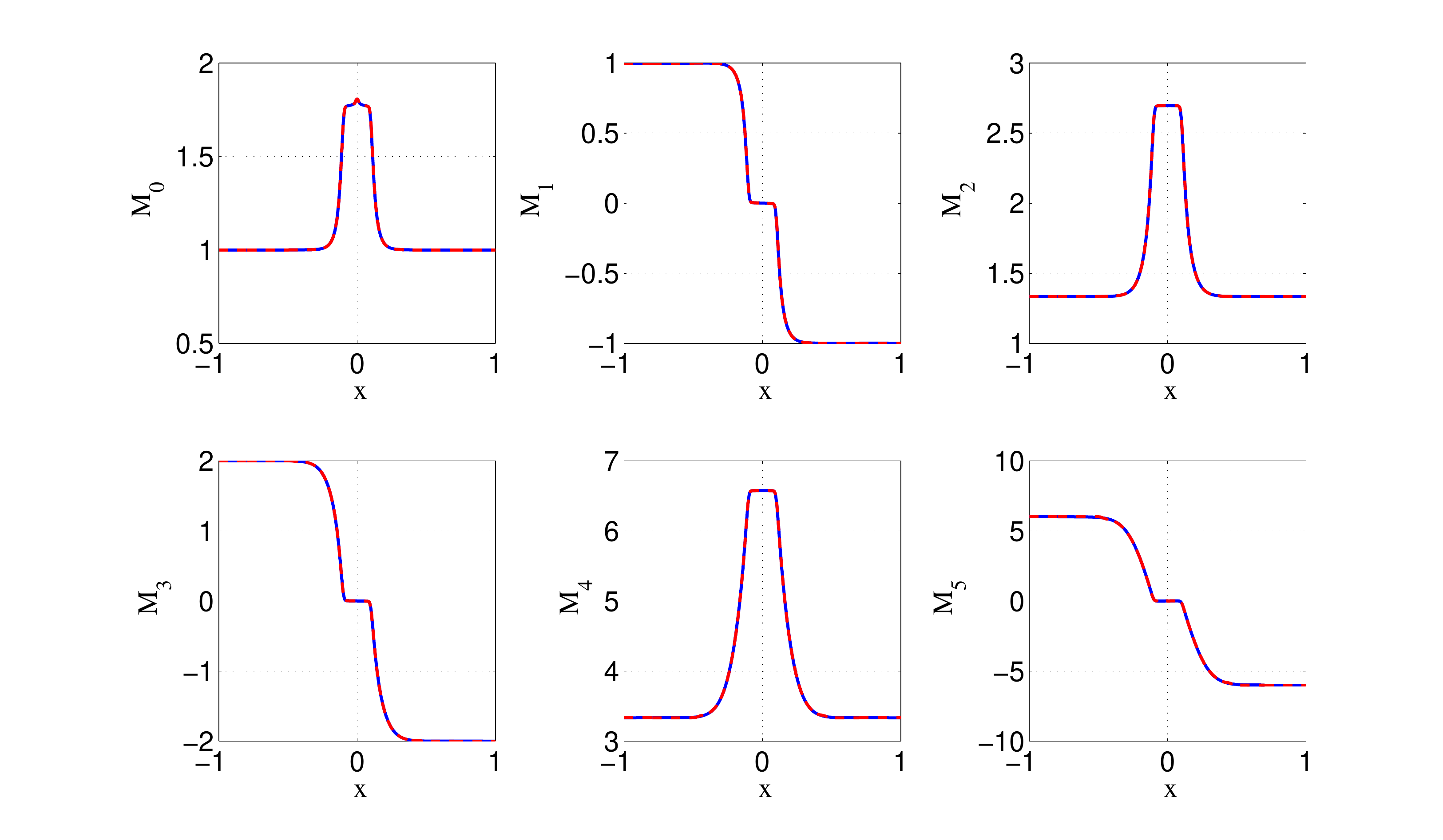}
\caption{Solution to 1-D Riemann problem at $t=0.5$. Five transported moments ($M_{0,...,4}$) and reconstructed moment ($M_5$). The dashed line is found with (\ref{eq:sigmar}) and the solid line with (\ref{eq:sigma}).\label{moments1D}}
\end{center}
\end{figure}

\begin{figure}[htb]
\begin{center}
\includegraphics[width = 1.0\textwidth,clip=true, trim= 40 0 40 0]{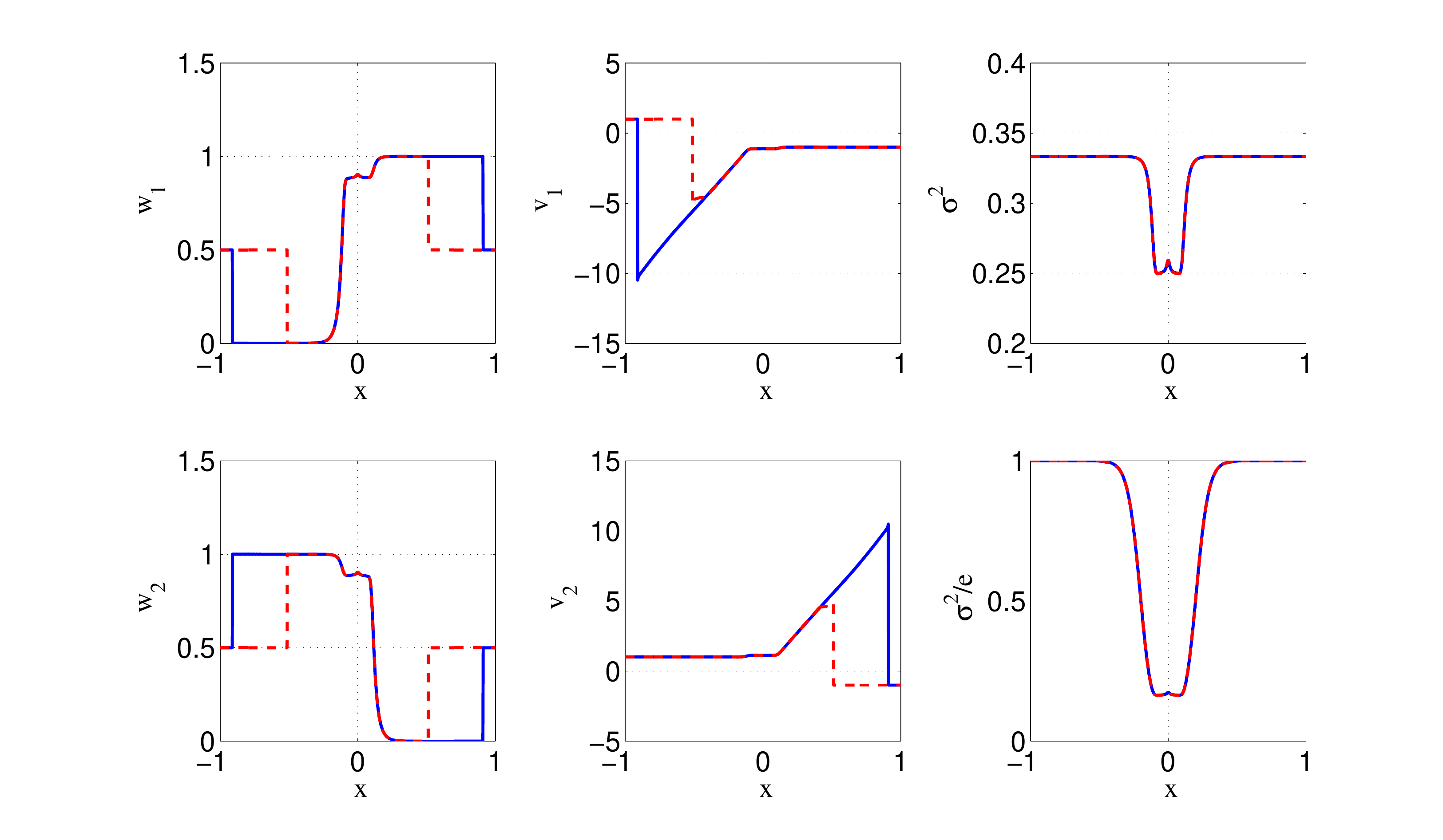}
\caption{Solution to 1-D Riemann problem at $t=0.5$. Left: weights. Center: abscissas. Top right: $\sigma^2$. Bottom right: Gaussian contribution to energy $\sigma^2/e$ (bottom right). The dashed line is found with (\ref{eq:sigmar}) and the solid line with (\ref{eq:sigma}).\label{quad1D}}
\end{center}
\end{figure}

\begin{figure}[htb]
\begin{center}
\includegraphics[width = 1.0\textwidth,clip=true, trim= 40 0 40 0]{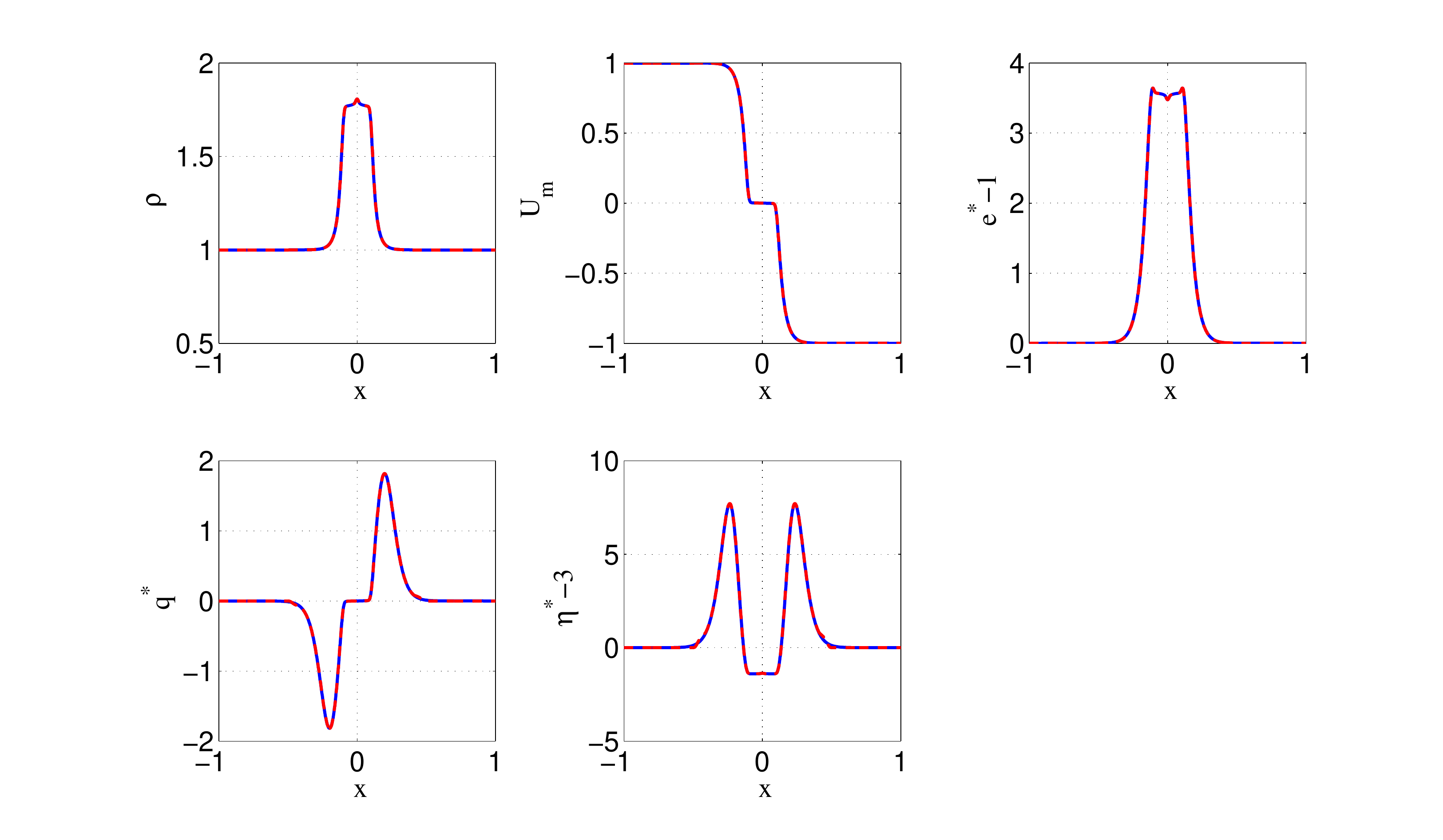}
\caption{Solution to 1-D Riemann problem at $t=0.5$. Top left: density $\rho$. Top center: mean velocity $U_m$. Top right: normalized energy $e^*-1$. Bottom left: skewness $q^*$. Bottom center: kurtosis $\eta^*-3$. The dashed line is found with (\ref{eq:sigmar}) and the solid line with (\ref{eq:sigma}).\label{centralmom1D}}
\end{center}
\end{figure}

Simulation results for the 1-D Riemann problem are presented in Figures~\ref{moments1D}--\ref{centralmom1D} at time $t=0.5$.  Note that due to the equilibrium initial conditions, only one velocity abscissa is used when $\sigma^2/e = 1$ (i.e., $v_1$) and the other ($v_2$) is set to zero automatically using the 1-D quadrature algorithm described in \S\ref{mia}. At $t=0.5$, one can observe from Figure~\ref{centralmom1D} that the equilibrium condition is still present on the left and right sides of the computational domain. In the center of the domain,  $\sigma^2/e \approx 0.2$, indicating that the overall distribution is composed of two Gaussian distributions with very little overlap. Also, note that unlike in a pure PTC problem where the velocity abscissas remain at their initial values (i.e., 1 and -1), in Figure~\ref{quad1D} the abscissas have their largest magnitudes just behind the ``shock'' in density at the edge of the equilibrium domain. This behavior is a direct result of the definition of the spatial fluxes in terms of the underlying two-node Gaussian-EQMOM distribution.  Indeed, the outer tails of the Gaussian distribution have higher velocity than the value at the peak density and thus penetrate faster into the equilibrium domain, resulting in a higher local flux velocity. The strong deviations from equilibrium are also clearly observed in the normalized energy, skewness, and kurtosis in Figure~\ref{centralmom1D}. Note, however, that using the velocity limiter in (\ref{eq:sigmar}) does not compromise the prediction of the moments in Figure~\ref{moments1D} because the high velocities correspond to negligibly small weights as can be seen in Figure~\ref{quad1D}.

Except at the edges of the equilibrium domain, we see from Figure~\ref{moments1D} that the transported moments and $\sigma^2$ are smoothly varying functions of $x$. More importantly, the singularities appearing in the solution do not belong to the class of $\delta$-shocks but to the less singular class of shocks encountered with hyperbolic systems of conservations laws, thus revealing a potentially well-behaved system. Moreover, due to the kinetic-based definition of the spatial fluxes, the moments are always realizable, and the moment-inversion algorithm always computes a well-defined quadrature from the updated moments. Overall, the two-node Gaussian-EQMOM reconstruction of the velocity distribution yields a robust numerical algorithm using a minimum number of moments. In comparison to the high-order delta function reconstruction described in \cite{fox09}, the two-node Gaussian-EQMOM provides a higher fidelity flux representation for a fixed number of transported moments. Moreover, because the moments of the Gaussian-EQMOM distribution can be computed to any desired order, the flux representation described in \S\ref{flux} can be systematically improved. This advantage becomes even more significant for 2-D and 3-D phase spaces where the number of transported moments needed for the delta-function reconstruction increases rapidly with the order of the moments \cite{fox09}.

\subsection{2-D particle-laden turbulent flow}\label{s:turbflow}

The proposed test case is a particle-laden gas phase represented by 2-D frozen homogeneous isotropic turbulence (HIT) generated with the ASPHODELE code of CORIA \cite{rev_cism07,reveillon07}, which solves the 2-D and 3-D low-Mach-number Navier-Stokes equations. The turbulence is generated following the Pope spectrum \cite{p00} with parameters $p_0=4$, $c_L=0.013$, $c_{\eta}=0.105$ and $\beta=5.2$. The particle phase is placed homogeneously in the domain at $t=0$ with the same velocity as the gas phase.  The droplet number density is uniform and the Stokes number range based on the Kolmogorov time scale ($t_\eta=0.3172$~s) is $\text{St}=\tau_p/t_\eta \in [1,20]$, which is large enough to observe particle trajectory crossings (PTC).  Predicting this type of flows is important, as it is expected to exhibit the main effects of a turbulence gas field on a disperse liquid spray, i.e., the preferential concentration of particles in low vorticity zones that greatly influences auto-ignition, which is of primary importance for turbulent combustion applications \cite{bouali2012}. This effect is highly size dependent, so the test case is particularly interesting for quantifying the accuracy of the proposed methods. 

\begin{figure}[htb]
\begin{center}
\includegraphics[width = 0.45\textwidth] {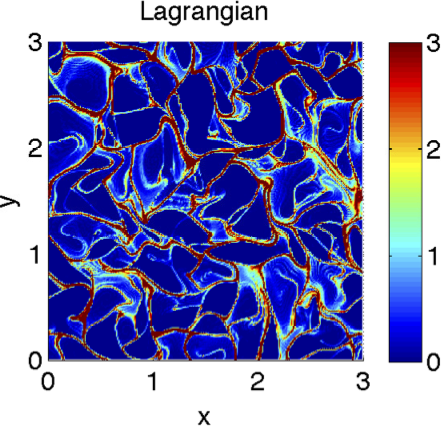}
\includegraphics[width = 0.45\textwidth] {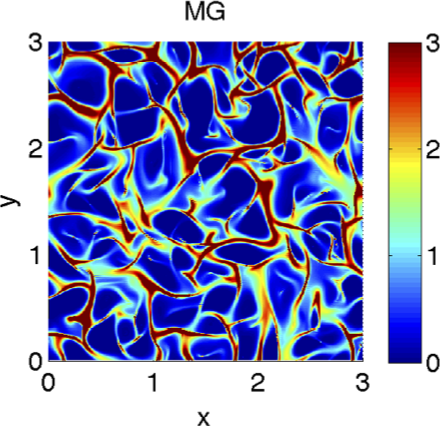}
\caption{Number density in frozen turbulence with $\text{St}=1$ at time $t=4$~s for Lagrangian (left) and Eulerian (right) simulations.}\label{fig:tf1}
\end{center}
\end{figure}

\begin{figure}[htb]
\begin{center}
\includegraphics[width = 0.45\textwidth] {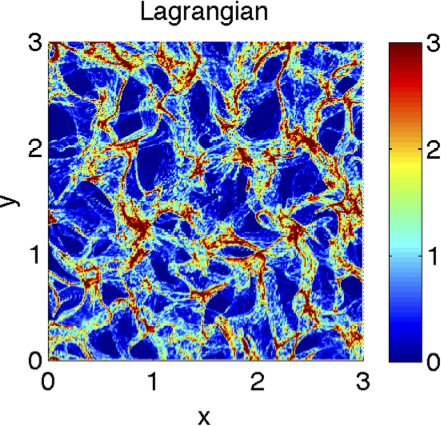}
\includegraphics[width = 0.45\textwidth] {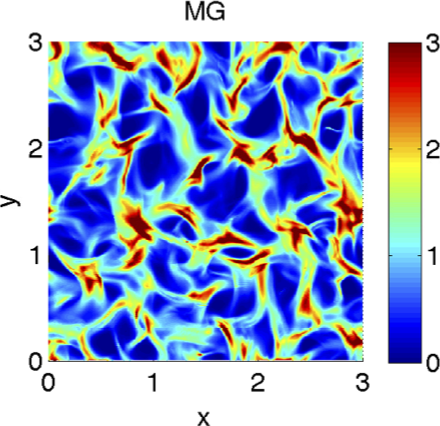}
\caption{Number density in frozen turbulence with $\text{St}=5$ at time $t=4$~s for Lagrangian (left) and Eulerian (right) simulations.}\label{fig:tf5}
\end{center}
\end{figure}

\begin{figure}[htb]
\begin{center}
\includegraphics[width = 0.45\textwidth] {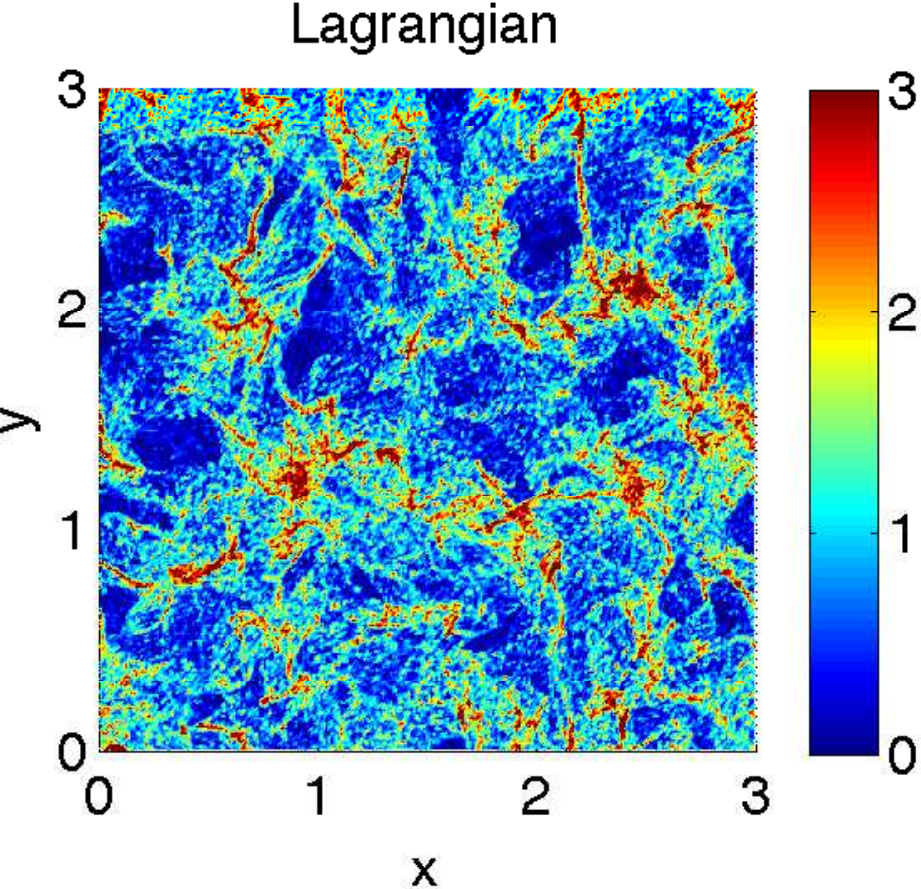}
\includegraphics[width = 0.45\textwidth] {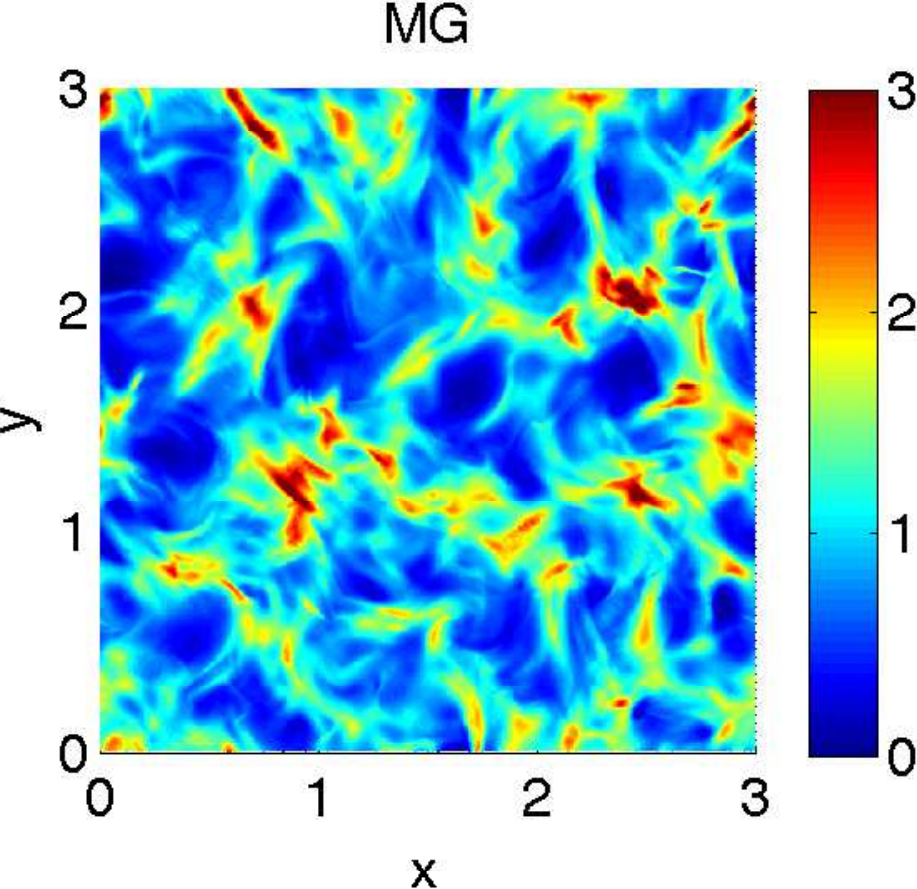}
\caption{Number density in frozen turbulence with $\text{St}=10$ at time $t=4$~s for Lagrangian (left) and Eulerian (right) simulations.}\label{fig:tf10}
\end{center}
\end{figure}

\begin{figure}[htb]
\begin{center}
\includegraphics[width = 0.45\textwidth] {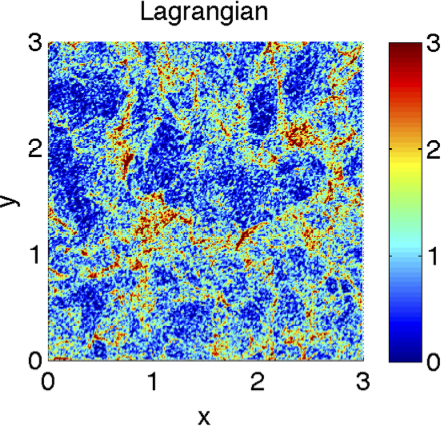}
\includegraphics[width = 0.45\textwidth] {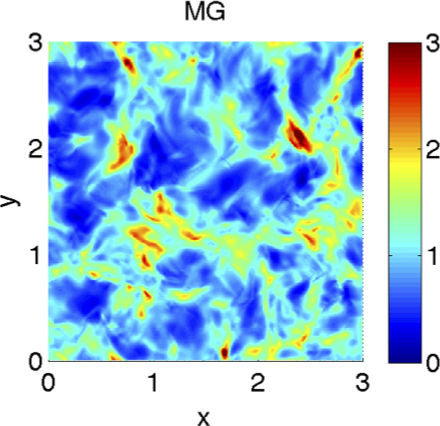}
\caption{Number density in frozen turbulence with $\text{St}=20$ at time $t=4$~s for Lagrangian (left) and Eulerian (right) simulations.}\label{fig:tf20}
\end{center}
\end{figure}

As a reference, Lagrangian computations using a first order in time scheme with time step $\Delta t=0.001$~s are performed using $10$ million particles, for which a satisfactory statistical convergence has been verified. The Eulerian projection of the Lagrangian particle statistics is done using a box filter at the grid-cell size.  Example comparisons of the number density field $M_{0,0}$ found with the Lagrangian and Eulerian (MG) simulations are shown for four different Stokes numbers in Figures~\ref{fig:tf1}--\ref{fig:tf20}, respectively, at $t=4$~s.  Note that at the time shown in the figures, the Eulerian fields have nearly reached a steady-state condition during which the moments change very little with time. This is possible because of the frozen (time-independent) nature of the turbulent gas phase.

At the smallest Stokes number shown in Figure~\ref{fig:tf1}, very little PTC is present.  However, very strong preferential concentration of the particles is evident.  In the absence of PTC, the Eulerian moment equations do not require a moment closure since the velocity distribution is monokinetic \cite{Sdc09_us} (i.e., the granular temperature is null).  Thus, the comparison between the Lagrangian and Eulerian number density fields can be used to judge the ability of the finite-volume scheme to capture the steep density gradients caused by the hypercompressibility of the particle phase.  From Figure~\ref{fig:tf1}, it is clearly evident that the Eulerian method is able to capture accurately the segregation of the particle phase, including regions in the domain where $M_{0,0}$ is null (i.e., vacuum zones), with very small numerical diffusion.

For $\text{St}=5$ shown in Figure~\ref{fig:tf5}, appreciable PTC is present and results in a decrease in the preferential concentration. Indeed, for this Stokes number, the velocity distribution has small (but nonzero) dispersion about the local mean velocity.  In the Lagrangian simulation, the full structure of the velocity dispersion is captured, while in the Eulerian model it is only partially captured due to the moment closure needed to close the system.  From Figure~\ref{fig:tf5}, we can observe that the Eulerian model captures much of the segregation structure seen in the Lagrangian number density field. Moreover, due to the hyperbolic nature of the moment closure, the Eulerian number density field is free from $\delta$-shocks caused by not resolving the PTC, which inevitably arise in weakly hyperbolic systems that do not adequately account for velocity dispersion \cite{kah11_cms}. 

For $\text{St}=10$ shown in Figure~\ref{fig:tf10}, the Lagrangian number density field becomes appreciably more uniform due to abundant PTC. Nevertheless, as clearly observed in the Eulerian field, near-vacuum zones still exist along side regions with preferential concentration of particles. Interestingly, due to the isotropic nature of the gas velocity field, the particle velocity dispersion becomes more isotropic with increasing Stokes number.  Under these conditions, the Eulerian simulation becomes less dependent on the choice of the moment closure (albeit as long as the chosen moment closure can represent velocity dispersion).  In fact, a full second-order moment closure based on a single Gaussian distribution \cite{vie2013cicp} (i.e., a degenerate case of Gaussian-ECQMOM) yields results very similar to Figure~\ref{fig:tf10} for homogeneous isotropic turbulence. However, for more complicated configurations with large-scale PTC (e.g., crossing particle jets), the full multi-Gaussian represention will be needed to capture the particle phase.

The last comparison for $\text{St}=20$ is shown in Figure~\ref{fig:tf20}.  Here, the Lagrangian number density field is free of vacuum zones since PTC is rampant.  Indeed, for this Stokes number the particles are accelerated slowly by the gas and, once they achieve a given velocity, they move in a manner that is weakly coupled to the gas velocity.  From the Eulerian field, we can observe that the moment closure again does a good job of reproducing the Lagrangian field. Furthermore, the Eulerian model has the advantage that the `stochastic noise' present in the Lagrangian field due to the finite sample of particles is absent in the moment model.  It would thus be possible to extract high-fidelity information concerning, for example, the spatial gradients of moments (e.g., energy spectra) that would be impossible to obtain from a Lagrangian simulation using a tractable number of particles.

From a mathematical perspective, the 2-D simulations presented above, carried out with the kinetic-based fluxes described in \S\ref{flux2}, provide evidence that the proposed multivariate Gaussian-EQMOM closure is robust for simulating particle-laden turbulent flows. In these simulations, no evidence of unrealizable moments was observed (although we have not proved that the kinetic-based fluxes in the form of (\ref{fluxfunc1a2D}) are realizable even for first-order schemes), nor was any evidence of weakly hyperbolic behavior when $V \neq 0$ (i.e., outside the conditions of Theorem~\ref{th_hyp2D}) observed.  

\section{Conclusions} \label{s:conclusion}

The multivariate Gaussian extended quadrature method of moments and the related moment-inversion algorithms appear to be a very promising approach for the direct-numerical simulation of particle-laden turbulent flows \cite{fox2012}.  The approach combines stability and a lower level of singularity compared to existing quadrature-based moment methods, see \cite{kah11_cms}, and is able to capture both particle trajectory crossing (PTC) caused by the free-transport term and the effects of turbulent agitation. It is noteworthy that the Gaussian-EQMOM naturally degenerates toward the correct velocity distribution with the associated spatial fluxes in both the PTC and dispersion limits. Moreover, by relying on the recent advances in CQMOM \cite{yuan2011}, the Gaussian-ECQMOM naturally adapts to the required number of nodes in even highly degenerate cases (e.g., in the absence of particles). As such, the Eulerian moment methods described in this work should offer an attractive alternative to Lagrangian particle tracking methods for simulating particle-laden flows.

\section*{Acknowledgments}

ROF was supported by a grant from the U.S.\ National Science Foundation (CCF-0830214). The research leading to the results reported in this work has received funding from the European Union Seventh Framework Programme (FP7/2007-2013) under grant agreement No.~246556.

\appendix

\section{Realizability of the first-order finite volume scheme with the kinetic-based flux}\label{s:realizability}

The finite volume scheme corresponding to the transport part of (\ref{modele_bipic}) can be written:
\begin{equation}\label{eq:schemeo1}
M_{k,j}^{n+1}  = M_{k,j}^{n} - \frac{\Delta t}{\Delta x_j} \left(F_{k,j+1/2}^n - F_{k,j-1/2}^n \right)
\end{equation}
where $M_{k,j}^{n}$ is an approximation at time $t^n$ of the mean value of the $k^{th}$-order moment on the cell $]x_{j-1/2},x_{j+1/2}[$ of size $\Delta x_j$.
Using the kinetic-based definition (\ref{fluxfunc}), the flux is defined by 
\begin{equation}\label{eq:fluxo1}
F_{k,j+1/2}^n = \int_{0}^{\infty} v^{k+1}f_j^n(v) \, \mathrm{d} v + \int_{-\infty}^0  v^{k+1} f_{j+1}^n(v) \, \mathrm{d} v, \quad k=0,\dots,4
\end{equation}
where $f_j^n$ is the two-node Gaussian-EQMOM reconstruction corresponding to the set of moments $(M_{k,j}^{n})_{k\in\{0,\dots,4\}}$ and given by the parameters $(\rho_{j,\alpha},v_{{j,\alpha}},\sigma_j)_{\alpha=1,2}$.

Let us denote by $(u_{\lambda,\alpha})_{\alpha=1,2,3}$ the three abscissas corresponding to the three-node quadrature of  the measure $\exp(-x^2)\mathds{1}_{]\lambda,+\infty[}(x) \mathrm{d} x$.
We then introduce the following conjecture, numerically checked for a large number of $\lambda\in\mathbb{R}$: 
\begin{itemize}
\item[{[C]}] $\forall \lambda\in\mathbb{R}, \quad \forall \alpha\in\{1,2,3\}, \qquad |u_{\lambda,\alpha}| \le \max_\alpha \left(u_{0,\alpha} \right) + \max\{0,\lambda\}$ .
\end{itemize}
Let us also remark that $\max_\alpha \left(u_{0,\alpha} \right)\le 1.8$.
The following proposition then gives a sufficient condition for the realizability of the scheme (\ref{eq:schemeo1}).
\begin{proposition}[\textrm{Realizability}] \label{prop_realizability}
Let us assume that the conjecture [C] is true. 
Then the scheme (\ref{eq:schemeo1}) with flux (\ref{eq:fluxo1}) is realizable if 
$$ \forall j, \quad \frac{\Delta t}{\Delta x_j}\max_\alpha (|v_{j,\alpha}|+1.8 \sigma_j\sqrt 2)\le 1.$$
\end{proposition}

\emph{Proof}. 
The scheme can be written, for $k\in\{0,\dots,4\}$:
$$
M_{k,j}^{n+1}  = \mathbf{I}_{j,k}
+\frac{\Delta t}{\Delta x_j}  \int_{-\infty}^0 v^{k}|v| f_{j+1}^n(v) \, \mathrm{d} v
+\frac{\Delta t}{\Delta x_j}  \int_{0}^{+\infty}  v^{k}|v| f_{j-1}^n(v) \, \mathrm{d} v,
$$
with
$$
\mathbf{I}_{j,k} =
M_{k,j}^{n} - \frac{\Delta t}{\Delta x_j} \int_{\mathbb{R}}v^{k} |v|  f_j^n(v) \, \mathrm{d} v  
=  \int_{\mathbb{R}}v^{k} \left(1 - |v| \frac{\Delta t}{\Delta x_j} \right)  f_j^n(v) \, \mathrm{d} v  
$$
It is sufficient to prove that $(\mathbf{I}_{j,k})_{k\in\{0,\dots,4\}}$ is a moment vector.
For that, let us use the reconstruction of $f_j^n$ and the formulas (\ref{fluxfunc1a},\ref{halfmom}):
$$
\begin{gathered}
\mathbf{I}_{j,k} =\sum_{\alpha=1}^{2} 
\frac{ \rho_{j,\alpha}}{\sqrt{\pi}} \left[
\int_{\frac{-v_{j,\alpha}}{\sqrt{2} \sigma_j}}^{\infty}  ( v_{j,\alpha} + \sqrt{2} \sigma_j s )^{k} 
\left(1-(v_{j,\alpha} + \sqrt{2} \sigma_j s)\frac{\Delta t}{\Delta x_j}\right)e^{-s^2}\, \mathrm{d} s \right.\\
\left. +\int^{\infty}_{\frac{v_{j,\alpha}}{\sqrt{2} \sigma_j}}  ( v_{j,\alpha} - \sqrt{2} \sigma_j s )^{k} 
\left(1+(v_{j,\alpha} - \sqrt{2} \sigma_j s)\frac{\Delta t}{\Delta x_j}\right)e^{-s^2} \, \mathrm{d} s
\right].
\end{gathered}
$$
Noticing that the integrals can be exactly computed by using the three-node quadratures of  the measures $\exp(-x^2)\mathds{1}_{]\frac{-v_{j,\alpha}}{\sqrt{2} \sigma_j},+\infty[}(x) \mathrm{d} x$ and $\exp(-x^2)\mathds{1}_{]\frac{v_{j,\alpha}}{\sqrt{2} \sigma_j},+\infty[}(x) \mathrm{d} x$, respectively, and using the conjecture concludes the proof.
\endproof

\bibliographystyle{siam}
\bibliography{./biblio}

\end{document}